\documentclass[letterpaper]{amsart}

\usepackage[T1]{fontenc}
\usepackage{color}
\usepackage{amsrefs}
\usepackage[utf8]{inputenc} 
\usepackage{comment}

\usepackage{mathtools}
\usepackage{colonequals}
\usepackage{amsmath}
\usepackage{amssymb}
\usepackage{tikz}
\usepackage{tikz-cd}
\usepackage{mathrsfs}
\usepackage{enumitem}
\usepackage[foot]{amsaddr}
\usepackage{mathtools}

\usepackage{mleftright}

\theoremstyle{plain}
\newtheorem{theorem}{Theorem}[section]
\newtheorem{corollary}[theorem]{Corollary}
\newtheorem{proposition}[theorem]{Proposition}
\newtheorem{lemma}[theorem]{Lemma}

\newtheorem*{theorem*}{Theorem}
\newtheorem*{lemma*}{Lemma}
\newtheorem*{prop*}{Proposition}
\newtheorem*{cor*}{Corollary}
\newtheorem*{conj*}{Conjecture}

\theoremstyle{definition}
\newtheorem*{definition*}{Definition}
\newtheorem{definition}[theorem]{Definition}

\theoremstyle{remark}
\newtheorem*{rem*}{Remark}
\newtheorem{rem}[theorem]{Remark}
\newtheorem*{example*}{Example}
\newtheorem{example}[theorem]{Example}

\numberwithin{theorem}{section}
\numberwithin{equation}{section} 
\numberwithin{figure}{section} 

\newcommand{\acknowledgement}{\subsection*{Acknowledgements}}
\newcommand{\Id}{\operatorname{Id}}
\newcommand{\R}{\mathbb{R}}
\newcommand{\Z}{\mathbb{Z}}

\newcommand{\N}{\mathbb{N}}

\newcommand{\catname}[1]{\mathbf{#1}}

\newcommand{\Hom}{\operatorname{\mathrm{Hom}}}
\newcommand{\Int}{\operatorname{\mathrm{Int}}}

\address{CONACYT-CIMAT, Centro de Investigaci\'on en Matem\'aticas, Calle Jalisco, Colon\'{i}a Valenciana, C.P. 36023 Guanajuato, GTO, Mexico}
\thanks{Research supported in part by European Research Council grant
FP7-ICT-318493-STREP, Cat\'{e}dras CONACYT / 1076, and Grant \#N62909-19-1-2134
from the US Office of Naval Research Global and the Southern Office of Aerospace Research and Development of the US Air Force Office of Scientific Research}

\title[A Unified Framework for Discrete and Continuous Homotopy]{\v Cech closure spaces:\\ A Unified Framework for Discrete and Continuous Homotopy}
\author{Antonio Rieser}
\date{}

\begin{document}

\begin{abstract}
Motivated by constructions in topological data analysis and algebraic combinatorics, we study homotopy theory on the category of \v Cech closure spaces $\catname{Cl}$, the category whose objects are sets endowed with a \v Cech closure operator and whose morphisms are the continuous maps between them. We introduce new classes of \v Cech closure structures on metric spaces, graphs, and simplicial complexes, and we show how each of these cases gives rise to an interesting homotopy theory. In particular, we show that there exists a natural family of \v Cech closure structures on metric spaces which produces a non-trivial homotopy theory for finite metric spaces, i.e. point clouds, the spaces of interest in topological data analysis. We then give a \v Cech closure structure to graphs and simplicial complexes which may be used to construct a new combinatorial (as opposed to topological) homotopy theory for each skeleton of those spaces. We further show that there is a Seifert-van Kampen theorem for closure spaces, a well-defined notion of persistent homotopy, and an associated interleaving distance. As an illustration of the difference with the topological setting, we calculate the fundamental group for the circle, `circular graphs', and the wedge of circles endowed with different closure structures. Finally, we produce a continuous map from the topological circle to `circular graphs' which, given the appropriate closure structures, induces an isomorphism on the fundamental groups.
\end{abstract}
\maketitle

\section{Introduction}
\label{sec:Introduction}
Homotopy theory has long been one of the primary tools used to study
topological spaces, and generalizations of the theory have had dramatic
implications in other areas as well, in particular in algebra and algebraic
geometry. There have recently been a number of attempts to extend the reach of
homotopy theory to more discrete geometrical objects, such as graphs
\cite{Babson_etal_2006}, directed graphs \cite{Grigoryan_et_al_2014}, and
simplicial complexes \cites{Barcelo_Kramer_Laubenbacher_2001,
Barcelo_Laubenbacher_2005,Barcelo_Capraro_White_2014}, and then to try to
characterize combinatorial properties of these objects in terms of their
discrete homotopy invariants. In parallel, a different approach to
discretization is developed in \cite{Plaut_Wilkins_2013}, in which the
homotopies themselves are discretized, and this is then used to show the
existence of certain relations in the fundamental group of geodesic spaces.
Ideas from algebraic topology are also being used to study spaces where the
natural topologies available don't capture the desired topological picture.
This occurs, for instance, when trying to infer information about the topology
of a manifold given a set of points sampled from it, a problem which has
motivated the development of topological data analysis
\cites{Carlsson_Zomorodian_2005, Carlsson_2009, Ghrist_2008}. (We remark in passing that, although it's true that there are many topologies on finite sets which have interesting homotopy groups, the neighborhoods in such topologies are typically unrelated to a metric on the set, which is undesirable when studying a set of points sampled from a metric space. For an extensive discussion of finite topological spaces, see \cite{Barmak_2011}.)

In this article, we develop homotopy theory in the category of closure spaces,
which, as we will show, allows for the application of homotopy theory in all of
the situations mentioned above, and which additionally reduces to the standard
theory for topological spaces. Although our approach is not necessarily
equivalent to the ones cited earlier, it nonetheless produces a unified
construction of homotopy theory in all of these different contexts, in addition
to defining a non-trivial homotopy theory for point clouds. Historically,
homotopies on closure spaces were first defined in
\cite{Demaria_Garbaccio_1984_2}, with the aim of developing an alternate
version of shape theory, which was pursued in \cites{Demaria_Garbaccio_1984, Demaria_Garbaccio_1985}. These homotopy groups were then applied to directed graphs in \cite{Demaria_1987}. A construction of homotopy groups on limit spaces (and therefore also closure spaces) is also given in \cite{Preuss_2002}, Section A.2, but was not developed further. To the best of our knowledge, this article is the first to apply the homotopy of closure structures to finite metric spaces and general skeleta of simplicial complexes, as well as the first to advance the general theory further since the works cited above. The plan of the paper is as follows. In Section 2, we give a formal introduction to closure spaces, and we collect the statements of results about them that we will need in what follows, Section 3 gives details of a specific family of closure operators of interest on metric spaces, in particular for topological data analysis on point clouds, and in Section 4 we study covering spaces, the fundamental group of the circle endowed with different closure structures, closure operators on graphs and skeleta of simplicial complexes, and give a definition of `persistent' homotopy and the related interleaving distance. In the remainder of this section, we introduce the core idea behind our approach in the context of point clouds.

Given a set of points sampled from a topological probability space $X$, it's
natural to ask whether the topological invariants of $X$ can be recovered from
the sample. Most current attempts at doing so, in particular persistent
homology and homological manifold learning, first assume that $X$ has a metric structure,
and then proceed by replacing the points with
balls of varying radii, effectively thickening the original set of
observations. In manifold learning \cites{Niyogi_Smale_Weinberger_2008,
Niyogi_Smale_Weinberger_2011, Chazal_et_al_2009}, the aim is to prove
that, under favorable conditions, i.e. with high probability given an
appropriate radius and a dense enough sample, the union of the balls centered
at the sample points is homotopy equivalent to the space from which the samples
were drawn, or else that some topological invariants of the original space may
be recovered. In persistent (\v Cech) homology \cites{Carlsson_2009, Ghrist_2008,Chazal_et_al_2016}, one attempts to recover topological information from a metric invariant built using a one-parameter family of unions of balls around the sample points. 

Instead of approximating the target space with auxiliary topological spaces built on our finite sample, however, our point of departure is to ask what sort of homotopy theory one may construct directly on finite sets of points, and then to construct (weak) homotopy equivalences between the sets of samples and the spaces from which they're sampled. While this might at first appear to lead to a trivial theory - and indeed it does, if one stays within the category $\catname{Top}$ - what we find is that, by changing the category, we are able to develop a homotopy theory which provides information about the global configuration of sample points without invoking an intermediate topological space. This is accomplished by `coarsening' the continuous maps rather than thickening the space, and, as we will see, this focus on the maps makes homotopy into a functor from the category of closure spaces to the category of groups. This, in particular, enables us to easily define a homomorphism between homotopy groups starting from a morphism in our category, something which remains difficult with other approaches, in particular for maps between metric spaces which do not preserve the metric.

We give a simple example to illustrate what is meant above by `coarsening' a continuous map. Consider the spaces $S^1 = [0,1]/(0\sim 1)$ and $X = \left\lbrace 0,\frac{1}{5},\dots,\frac{4}{5}\right\rbrace \subset S^1$. We would like to know what properties of $S^1$ we can recover from properties of $X$. From a homotopical point of view, we would like to find a category in which there is an equivalence between $S^1$ and $X$ that we may invert in the homotopy category, which means, at the very least, that we need to have non-trivial morphisms from $S^1$ to $X$, even though, topologically, the only continuous maps are the constant maps. In closure spaces, however, we have a new class of maps to consider. In Section \ref{sec:ScMet}, we show that, on a family of closure structures on metric spaces, there is an $\epsilon$-$\delta$ version of continuity for these closure spaces, which we describe here. For a pair of non-negative real numbers $(q,r)$, we say that a map $f:S^1 \to X$ is $(q,r)$-continuous at a point $x$ iff $\forall\epsilon>0, \exists \, \delta(x)>0$ such that $\Vert x - y \Vert < q + \delta$ implies $\Vert f(x)- f(y)\Vert < r + \epsilon$. According to this defintion, the `nearest neighbor map' $f:S^1 \to X$ given by
\begin{equation}
\label{eq:Nearest neighbor map}
f(x) = \begin{cases}
0, & \text{for }x \in \left(\frac{9}{10},1\right) \cup \left[0,\frac{1}{10}\right]\\
\frac{1}{5}, & \text{for }x\in \left(\frac{1}{10},\frac{3}{10}\right]\\
\frac{2}{5}, &\text{for } x \in \left(\frac{3}{10},\frac{5}{10}\right]\\
\frac{3}{5}, & \text{for }x \in \left(\frac{5}{10},\frac{7}{10}\right]\\
\frac{4}{5}, & \text{for }x \in \left(\frac{7}{10}, \frac{9}{10} \right]
\end{cases}
\end{equation}
is $(q,r)$-continuous for $r \geq \frac{1}{5}$ and $0 \leq q \leq r-\frac{1}{5}$, although it is clearly discontinuous topologically. To see this continuity, note that, seen as a subset of $S^1$, each point in $X$ is at most a distance of $\frac{1}{10}$ from any given point in $S^1$. We therefore have
\begin{align*}
|f(x) - f(y)| &\leq |f(x) - x| + |x - y| + |y - f(y)| \\
&< \frac{1}{10} + q + \delta + \frac{1}{10} \leq r + \epsilon.
\end{align*} 
for any $(q,r)$ with $q + \frac{1}{5} \leq r$, where we let $\delta=\epsilon$. As we will see in Sections \ref{sec:ScMet} and \ref{sec:HomConv}, we may also construct homotopies and homotopy groups so that, for a certain range of $(q,r)$, both spaces have the same fundamental group as the (topological) circle. This example illustrates two important points. The first is that topologically continuous maps are, in some sense, too rigid to be useful for a homotopy theory of point clouds, and the second is that they are not the only choice. 

 We make three additional observations about the modified notion of continuity above. First, note that for the pair $\left(q=0,r=\frac{1}{10}\right)$, the map $f$ is not $\left(0,\frac{1}{10}\right)$-continuous. In particular, this illustrates how modifying continuity in this way allows (topological) discontinuities at any point, but only if the jumps at (topologically) discontinuous points are controlled, the allowable size of the (topological) discontinuity being given by the pair $(q,r)$.

Second, we observe that topological continuity does not imply $(q,r)$-continuity for all $q$ and $r$. Consider, for instance, multiplication by $2$ on the real line, i.e. $f:\R \to \R$, $f(x)\coloneqq 2x$. While clearly continuous topologically, $f$ is not $(q,r)$-continuous for any pair $(q,r)$ where $q>\frac{1}{2}r>0$. To see this, first observe that any interval $I_{x,r,\epsilon} \coloneqq (x-r-\epsilon,x+r+\epsilon)$ must necessarily be the image of the interval $\left(\frac{1}{2}(x-r-\epsilon),\frac{1}{2}(x+r+\epsilon)\right)$. It's therefore clear that, for $q>\frac{1}{2}r$ and sufficiently small $\epsilon > 0$, there is no $\delta>0$ for which the interval $f\left(\frac{x}{2}-q-\delta,\frac{x}{2}+q+\delta\right) \subseteq I_{x,r,\epsilon}$, and so $f$ is not $(q,r)$-contiuous at the point $\frac{x}{2}$. This illustrates that, while local discontinuities are allowed, a $(q,r)$-continuous function must be rather uniformly controlled at the scale which determines $(q,r)$-continuity. While at first perhaps unsettling, we will see that, for point clouds, this rigidity has the desirable effect of keeping nearby points near each other after the application of a $(q,r)$-continuous function.

Finally, we remark that $(q,r)$-continuous maps between metric spaces need not be coarse maps in the sense of coarse geometry. That is, at large enough scales, a $(q,r)$-continuous map may send points which are initially a finite distance apart arbitrarily far away from each other, so long as the local $(q,r)$-continuity condition is satisfied. As an example, take the set of points
\[ X = \{(-1,i)\}_{i\in \N \cup {0}} \cup \{(0,0)\} \cup \{(1,j)\}_{j\in \N \cup {0}}.\]
Define the function $F:X \to \R^2$ by $F(a,b) \coloneqq (a\cdot b,0)$. This
function is $(1,1)$-continuous, but it is not a coarse map, since there is no
uniform bound on the distance $d(F(a,b),F(c,d))$. Together, these examples
illustrate the fundamental features of the maps between the closure structures
on metric spaces that will be our main objects of interest: they are wild at small scales, rigid at medium scales, and flexible at large scales.

\section{The category $\catname{Cl}$}
\label{sec:Closure spaces}

A \emph{closure structure} on a space $X$ collects all of the information about the neighborhoods of each point in $X$. It is weaker than a topology, but, as we will see, still allows for the construction of a rich homotopy theory which extends classical homotopy on the subcategory $\catname{Top}$. Among other things, this will allow us to construct weak homotopy equivalences between spaces which are topologically very different, but where there exist closure structures with similar characteristics, giving rise to isomorphic homotopy groups. In this section, we define closure spaces, and we include the statements of a number of results which we will use in this paper, both as an introduction to the subject - which, while classical, is not broadly known - as well as for reference.

In the examples and computations in this article, we will be mainly concerned with a natural family of closure structures induced by a metric, but the theory holds unchanged in the general setting. For additional general results on closure spaces and for the proofs of the results given here, we refer the reader to the results on closure spaces in the book \cite{Cech_1966}.

\subsection{\v Cech closure spaces}

\begin{definition}
\label{defn:Closure structure}
Let $X$ be a set. A \emph{\v {C}ech closure operator} on $X$ is a map $c:\mathcal{P}(X) \to \mathcal{P}(X)$ which satisfies
\begin{enumerate}[itemindent=*, leftmargin=*]
\item $c(\emptyset) = \emptyset$
\item $A \subseteq c(A)$
\item $c(A \cup B) = c(A) \cup c(B)$.
\end{enumerate} 
A pair $(X,c)$ is called a \emph{\v Cech closure space}, or simply a \emph{closure space}, and for a set $A \subset X$, we call $c(A)$ the \emph{closure} of $A$. If $A = \{x\}$, we will write $c(\{x\})$ as $c(x)$. We sometimes say that the map $c$ gives a \emph{\v {C}ech closure structure} on $X$.

Given two \v Cech closure operators $c_1$ and $c_2$ on the same space $X$, we say that \emph{$c_1$ is finer than $c_2$}, and \emph{$c_2$ is coarser than $c_1$}, iff, for each $A \subset X$, $c_1(A) \subset c_2(A)$, and we write $c_1 < c_2$ to denote this relation. 
\end{definition}

\begin{rem}
Note that the definition above immediately implies that the closure operator is monotone, i.e. $A\subset B \implies c(A) \subset c(B)$. \end{rem}

We will make use of the following special closure structures in the sections that follow.

\begin{definition}
	Let $X$ be a set. If a closure structure $c$ on $X$ satisfies $c(x) = x$ for every point $x \in X$, then we say that $c$ is the \emph{discrete closure structure} on $X$, or simply \emph{discrete}. If a closure structure $x$ on $X$ satisfies $c(A) = X$ for every nonempty $A \subset X$, then we say the $c$ is the \emph{indiscrete closure structure} on $X$, or simply \emph{indiscrete}.
\end{definition}

We may also construct closure spaces via interior operators, defined as follows.

\begin{definition}
Let $X$ be a set. An \emph{interior operator} on $X$ is a map $\Int:\mathcal{P}(X) \to \mathcal{P}(X)$ which satisfies
\begin{enumerate}
\item $\Int(X) = X$
\item $\Int(A) \subseteq A$
\item $\Int(A \cap B) = \Int(A) \cap \Int(B)$.
\end{enumerate}
\end{definition}

\begin{definition}
	\label{def:I from c and c from i}
Given a closure space $(X,c)$ we define the operator $i_c:\mathcal{P}(X) \to \mathcal{P}(X)$ to be
\begin{equation*}
i_c(A) = X - c(X - A).
\end{equation*}
Similarly, given a set $X$ with an interior operator $i$, we define an operator $c_i:\mathcal{P}(X) \to \mathcal{P}(X)$ by
\begin{equation*}
c_i(A) = X - i(X - A)
\end{equation*} 
\end{definition}

\begin{proposition}
For any closure space $(X,c)$, $i_c$ is an interior operator on $X$, and $c_{(i_c)} = c$. Similarly, for any interior operator $i$ on a set $X$, $c_i$ is a \v Cech closure operator on $X$, and $i_{(c_i)} = i$.
\end{proposition}
\begin{proof}
Given a closure operator $c$ on $X$, it follows immediately from the axioms for closure operators and de Morgan's laws that $i_c$ is an interior operator. In the same way, given an interior $i$ operator on $X$, $c_i$ is a closure operator by the axioms for interior operators and de Morgan's laws. Finally, for any $A \subseteq X$, we have
\begin{align*}
c_{(i_c)}(A) &= X - i_c(X - A) = X - (X - c(X - (X - A))) = c(A), \text{ and }\\
i_{(c_i)}(A) &= X - c_i(X - A) = X - (X - i(X - (X - A))) = i(A).
\end{align*}
\end{proof}

\begin{definition}
	Let $(X,c)$ be a closure space. We say that a subset $A \subset X$ is \emph{closed} if $c(A)=A$, and that $A \subset X$ is \emph{open} if $X - A$ is closed. If $\emptyset \neq A \subsetneq X$ is open, then we say that $A$ is a \emph{proper open subset} of $X$.
\end{definition}

While open and closed sets exist in closure spaces, in the following they take a secondary role to the neighborhoods of a set, defined below, which are not necessarily open or closed.

\begin{definition}
Let $(X,c)$ be a closure space. We say that a set $U \subseteq X$ is a
\emph{neighborhood} of a set $A \subseteq X$ if $A \subseteq X - c(X-U)$. The
\emph{neighborhood system} of a set $A$ is the collection of neighborhoods of $A$. 
\end{definition}

Neighborhoods and open sets are related by the following theorem.

\begin{theorem}[\cite{Cech_1966}, Theorem 14.B.2]
\label{thm:Open subset characterization}
Let $(X,c)$ be a closure space. A subset $U\subset X$ is a neighborhood of a subset $A \subset X$ iff $U$ is a neighborhood of each point of $A$. A subset $U \subset X$ is open in $(X,c)$ iff it is a neighborhood of all of its points, or, equivalently, if it is a neighborhood of itself.
\end{theorem}

\begin{definition}
	Let $(X,c)$ be a closure space. A \emph{base of the neighborhood system of $A \subset X$} is a collection $\mathcal{B}$ of subsets of X such that each set $B \in \mathcal{B}$ is a neighborhood of $A$, and each neighborhood of $A$ contains a set in $\mathcal{B}$. A \emph{subbase of the neighborhood system of $Y \subset X$} is a collection $\mathcal{C}$ of subsets of $X$ such that the collection of all finite intersections of elements of $\mathcal{C}$ is a base of the neighborhood system of $Y$. When $A$ contains only a single point $x \in X$, we will sometimes use the term \emph{local base (local subbase) at $x$} to refer to the base (subbase) of the neighborhood system of $\{x\}$.
\end{definition}

The next theorems show that we may also construct a closure structure on a set $X$ by specifying the local bases at each point. We first note, however, that any local base satisfies the following

\begin{proposition}[\cite{Cech_1966}, 14.B.5]
\label{prop:Local base properties}
Let $(X,c)$ be a closure space. If $\mathcal{U}(x)$ is a local base at a point $x$, then the following conditions are satisfied:
\begin{enumerate}
\item \label{item:Local base point 1}$\mathcal{U}(x) \neq \emptyset$
\item For each $U \in \mathcal{U}(x)$, $x \in U$
\item \label{item:Local base point 3}For each $U_1, U_2 \in \mathcal{U}(x)$, there exists a $U\in \mathcal{U}(x)$ such that $U \subset U_1 \cap U_2$.
\end{enumerate}
\end{proposition}

We now see how to construct closure structures on a set $X$ by designating a special collection of sets for each point $x \in X$.

\begin{theorem}[\cite{Cech_1966}, Theorem 14.B.10]
\label{thm:Closure from collection of sets}
For each point $x \in X$, let $\mathcal{U}(x)$ be a collection of sets that satisfies conditions \ref{item:Local base point 1}-\ref{item:Local base point 3} in Proposition \ref{prop:Local base properties}. Then there exists a unique closure structure $c$ on $X$ such that, for each $x \in X$, $\mathcal{U}(x)$ is a local base at $x$ in $(X,c)$.
\end{theorem}

\begin{corollary}[\cite{Cech_1966}, Corollary 14.B.11]
\label{cor:Closure with prescribed base and subbase}
\begin{enumerate}[wide, labelindent = 0pt]
Let $X$ be a set, and for any non-empty collection of subsets $\mathcal{U}$ of $X$, let $\bigcap \mathcal{U}$ denote the intersection of the elements of $\mathcal{U}$, i.e.
\begin{equation*}
\bigcap \mathcal{U} \coloneqq \bigcap_{U \in \mathcal{U}} U
\end{equation*}
\item For each element $x \in X$, let $\mathcal{U}(x)$ be a filter on $X$ such that $x \in \bigcap \mathcal{U}(x)$. Then there exists a unique closure structure $c$ on $X$ such that $\mathcal{U}(x)$ is the neighborhood system at $x$ in $(X,c)$ for every $x \in X$.
\item \label{item:Closure with prescribed subbase} For each element $x \in X$, let $\mathcal{V}(x)$ be a non-empty family of subsets of $X$ such that $x \in \bigcap \mathcal{V}(x)$. Then there exists a unique closure structure $c$ on $X$ such that $\mathcal{V}(x)$ is a local subbase at $x$ in $(X,c)$ for every $x \in X$.
\end{enumerate}
\end{corollary}

The next theorem and corollary show the connection between the closure of a set and neighborhoods of points, generalizing familiar facts from topological spaces.

\begin{theorem}[\cite{Cech_1966}, Theorem 14.B.6]
	\label{thm:Point contained in closure if neighborhood intersects}
	Let $(X,c)$ be a closure space. Given a set $A \in X$, a point $x\in X$ is contained in $c(A)$ iff every neighborhood $U$ of $x$ intersects $A$ non-trivially.
\end{theorem}
	
\begin{corollary}[\cite{Cech_1966}, Corollary 14.B.7]
	\label{cor:Point contained in closure if local base intersects}
	Let $(X,c)$ be a closure space. If $\mathcal{U}$ is a local base at a point $x \in X$, then for a subset $A\subset X$, $x \in c(A)$ iff, for each $U \in \mathcal{U}$, $U \cap A \neq \emptyset$.
\end{corollary}

\begin{rem}
	Note that, by Theorem \ref{thm:Point contained in closure if neighborhood intersects} and Corollary \ref{cor:Point contained in closure if local base intersects}, the closure operator $c$ constructed in Theorem \ref{thm:Closure from collection of sets} and Corollary \ref{cor:Closure with prescribed base and subbase} must be of the form
	\begin{equation*}
	c(A) = \{x \in X \mid \forall U \in \mathcal{N}(x), U\cap A \neq \emptyset\},
	\end{equation*}
	where $\mathcal{N}(x)$ is the neighborhood system of $x$.
\end{rem}

\begin{example}
Let $X = \{0,1,2,3\}$ and $c(x) = \{(x - 1) ,x, (x+1) \}$, where the integers are interpreted modulo $4$. Then $\{(x - 1),x, (x+1)\}$ is a neighborhood of $\{x\}$ for every $x \in X$.
\end{example}

\begin{example}
Every topological space $X$ is a closure space with the closure operator
defined by $c(A) = \bar{A}$. Note that for the closure operator on topological
spaces, $c(c(A)) = c(A)$. Furthermore, if $(X,c)$ is a closure space with $c^2
= c$, then the collection of sets $\mathcal{U} \coloneqq \{U \subset X \mid X-U
\text{ is closed, i.e. } c(X-U)=X-U\}$ forms a topology. To see this, first
note that $\emptyset$ and $X$ are in $\mathcal{U}$, since $c(\emptyset) =
\emptyset$, and $X \subset c(X), c(X) \subset X \implies c(X)=X$. Second, for
an arbitrary collection of sets in $\mathcal{U}$, say $\mathcal{U}'\coloneqq
\{U_{\lambda} \subset X \mid \lambda \in \Lambda\}$, we have $X - \cup_{\lambda
\in \Lambda} U_\lambda \subseteq c(X - \cup_{\lambda \in \Lambda}U_{\lambda})$.
Conversely, since closure operators are monotone, $c(\cap_{\lambda\in
\Lambda}(X - U_\lambda)) \subseteq c(X-U_{\lambda}) = X - U_{\lambda}$ for each $U_\lambda$, and therefore 
\begin{equation*}
c(X - \cup_{\lambda \in \Lambda}U_{\lambda}) = c(\cap_{\lambda\in \Lambda} (X - U_\lambda)) \subseteq \cap_{\lambda \in \Lambda} X-U_\lambda = X - \cup_{\lambda\in \Lambda} U_\lambda. 
\end{equation*} 
It follows that $\cup_{\lambda \in \Lambda} U_\lambda \in \mathcal{U}$. Finally, for a finite intersection of sets $\{U_i\}_{i=1}^n$ in $\mathcal{U}$, $\cap_{i=1}^n U_i$, we have
\begin{equation*}
c(X - \cap_{i=1}^n U_i) = c(\cup_{i=1}^n X-U_n) = \cup_{i=1}^n c(X - U_i) =  \cup_{i=1}^n X - U_i = X - \cap_{i=1}^n U_i.
\end{equation*}
Therefore, $\cap_{i=1}^n U_i$ is in $\mathcal{U}$, so $\mathcal{U}$ forms the
open sets of a topology whose closed sets are the fixed points of the operator
$c$. A closure operator satisfying $c^2 = c$ is called a \emph{topological} or \emph{Kuratowski closure operator}.
\end{example}

Kuratowski closure operators are characterized by the following

\begin{theorem}[\cite{Cech_1966}, Theorem 15.A.2]
	Each of the following conditions is necessary and sufficient for a closure space $(X,c)$ to be topological:
	\begin{enumerate}
		\item The closure of each subset of $(X,c)$ is closed in $(X,c)$
		\item The interior of each subset is of $(X,c)$ is open in $(X,c)$
		\item \label{thm:Open sets are local base in top spaces} For each $x \in X$, the collection of all open neighborhoods of $x$ is a local base at $x$.
		\item For each $x \in X$, if $U$ is a neighborhood of $x$, then there exists a neighborhood $V$ of $x$ such that $U$ is a neighborhood of each point of $V$, i.e. every neighborhood of any point $x \in X$ is a neighborhood of a neighborhood of $x$.
	\end{enumerate}
\end{theorem}

\begin{rem}
A pair $(X,c_X)$ is sometimes called a \emph{pretopology} in the literature. However, since these pairs, and not topologies, are our main objects of study, we have elected to revert to the older nomenclature used in \cite{Cech_1966}, which we believe does not semantically relegate these spaces to secondary, or preparatory, status. This convention also has the advantage of making \v {C}ech closure structures terminologically distinct from pretopologies in the sense of Grothendieck, which are different objects altogether. In this article, a closure space will always refer to a space with a closure operator that satisfies the axioms in Definition \ref{defn:Closure structure}, although the reader should be warned that there is some variation in the literature. In particular, except when explicitly stated, we will not require that $c^2 = c$, which, combined with the axioms in Definition \ref{defn:Closure structure}, would make the closure operator into what is known as a \emph{Kuratowski closure operator}. As Kuratowski closure operators induce a topology whose closed sets are the closed sets given by the operator, we will refer to spaces with Kuratowski closure operators simply as topological spaces.
\end{rem}

\subsection{Covering systems}
 \v {C}ech closure structures are most interesting when, as in the case of many
 examples constructed in Section \ref{sec:ScMet}, there are too few open sets
 to form a rich topology. In such cases, we will find that the neighborhood is,
 in fact, the more important object, and that there are many more neighborhoods
 than open sets. We now introduce covering systems in the context of closure spaces, one of the most important manifestations of the difference between neighborhoods and open sets.

\begin{definition}
Let $\mathcal{C}=\{U_\alpha\}_{\alpha\in I}$ be a family of subsets of a
closure space $(X,c)$. We say that $\mathcal{C}$ is a \emph{cover} of $(X,c_X)$
if $\cup_{U \in \mathcal{C}} U = X$, i.e. if every point $x \in X$ is contained in
some $U \in \mathcal{C}$. We say that $\mathcal{C}$ is an \emph{interior cover}
of $(X,c_X)$ if $\cup_{U \in \mathcal{C}} i_c(U) = X$, i.e. if every point $ x \in X$ has a neighborhood in $\mathcal{C}$.
\end{definition}

\begin{example}
\begin{enumerate}[wide]
\item Any cover of a closure space $(X,c)$ by open sets is an interior cover, which we call an \emph{open cover} of $(X,c)$.
\item Let $(V,E)$ be the graph defined by \[V = \{0,1,2,3\}, E = \{(k,k+1) \mid
    k \in V\},\] where the integers are to be understood modulo $4$. Then the
	family \[ \{ \{k-1, k, k+1\} \mid k \in V \}\] is an interior cover of
	the closure space
	$(V,c_E)$, where $c_E$ is defined by
	\begin{align*}
		&c_E(v) \coloneqq \{ w \mid w = v \text{ or } (v,w)\}, \text{ for
			every vertex } v \in V,\\
		&c_E(A) \coloneqq \cup_{v \in A} c_E(v) \text{ for any subset } A \subset V.
	\end{align*}
\item For the graph $(V,E)$ above, the family of sets $\mathcal{C} = \{ \{i-1,
    i, i+1\} \mid i \in \{0,2\} \}$  is a cover of the set $V$, but $\mathcal{C}$ is not an interior cover of the closure space $(V,c_E)$, since there are no neighborhoods of the vertices $v=1$ or $v=3$ in $\mathcal{C}$. 
\end{enumerate}
\end{example}

Compactness in closure spaces now takes the following form.

\begin{definition}
\label{def:Compactness}
We say that a closure space $(X,c)$ is \emph{compact} iff every interior cover
of $X$ has a finite subcover.
\end{definition}

\begin{rem*}
    We note that the finite subcover in the Definition \ref{def:Compactness} above may not itself be an
    interior cover. We also remark that, while Definition 41.A.3 in \cite{Cech_1966} defines compactness on closure spaces
    in terms of filters, the above definition is equivalent by
    Theorem 41.A.9 in \cite{Cech_1966}, and this one is the more useful version for
    the purposes of this article. 
\end{rem*}

\subsection{Continuous functions}

Continuity for maps between closure spaces is defined as follows. 

\begin{definition}
Let $(X,c_X)$ and $(Y,c_Y)$ be closure spaces. We say that a map $f:X\to Y$ is \emph{continuous at x} iff, for any subset $U\subseteq X$, $x\in c(U) \implies f(x) \in c(f(U))$. We say that a function $f$ is \emph{continuous} iff $f$ is continuous at every point $x$ of $X$. Equivalently, a function $f$ is continuous iff for every set $U\subseteq X$, $f(c(U)) \subseteq c(f(U))$.
\end{definition}

We present in this subsection several important basic results on continuous
maps between closure spaces. The proofs and additional results may be found in
\cite{Cech_1966}, Section 16.A. We begin with a result showing that the composition of
continuous functions is continuous.

\begin{proposition}[\cite{Cech_1966}, 16.A.3]
	\label{prop:Composition of continuous is continuous}
	Let $(X,c_X)$, $(Y,c_Y)$, and $(Z,c_Z)$ be closure spaces. Suppose $f:(X,c_X) \to (Y,c_Y)$ is continuous at $x \in X$, and $g:(Y,c_Y) \to (Z,c_Y)$ is continuous at $f(x)$, then the composition $g \circ f:(X,c_X) \to (Z,c_Z)$ is continuous at $x$. In particular, the composition of two continuous maps is continuous.
\end{proposition}

The next proposition indicates when the identity is continuous, as a map on the same space with two different closure structures.

\begin{proposition}[\cite{Cech_1966}, 16.A.2]
\label{prop:Identity cont}
Let $c$ and $c'$ be closure operations on a set $X$. Then $c < c'$ (i.e. $c$ is finer than $c'$) iff the identity map $id:(X,c) \to (X,c')$ is continuous.
\end{proposition}

The next result gives a characterization of continuity in terms of neighborhoods.

\begin{theorem}[\cite{Cech_1966}, Theorem 16.A.4 and Corollary 16.A.5]
\label{thm:Neighborhood continuity}
Let $(X,c)$, and $(Y,c')$ be closure spaces. A map $f:(X,c) \to (Y,c')$ is continuous at $x \in X$ iff, for every neighborhood $V \subset Y$ of $f(x)$, the inverse image $f^{-1}(V) \subset X$ is a neighborhood of $x$. Equivalently, $f$ is continuous at $x$ iff, for each neighborhood $V \subset Y$ of $f(x)$, there exists a neighborhood $U \subset X$ of $x$ such that $f(U) \subset V$. The function $f$ is continuous iff the above is satisfied for every point $x \in X$.
\end{theorem} 

The above result further specializes to open and closed sets.

\begin{corollary}[\cite{Cech_1966}, 16.A.6] If $f:(X,c) \to (Y,c')$ is a continuous mapping between \v Cech closure spaces, then for any open (closed) set $V \subset Y$, the preimage $f^{-1}(V) \subset X$ is open (closed). 
\end{corollary}

In particular, this implies

\begin{corollary}
Let $(X,c)$ be a closure space whose only open sets are $\emptyset$ and $X$, and suppose that $(Y,\tau)$ is a $T_0$ topological space. Then every continuous map $f:(X,c) \to (Y,\tau)$ is constant.
\end{corollary}
\begin{proof} Suppose $f$ is not constant. Then there are two points $x,y \in
	X$ such that $f(x) \neq f(y)$.  Since $Y$ is $T_0$, every pair of
	points $x,y \in Y$, at least one of them has an open neighborhood $U$
	which does not contain the other point. Without loss of generality, let
	$V \subset Y$ be a neighborhood of $f(x)$ with $f(y) \notin V$. Then
	$f^{-1}(V) = \emptyset$ or $X$, which is a contradiction, since, on the
	one hand, $x \in f^{-1}(V)$, so $f^{-1}(V) \neq \emptyset$, and, on the
	other hand, $y \notin f^{-1}(V)$, so $f^{-1}(V) \neq X$. Therefore, $f$ is constant.   
\end{proof}

\begin{rem}
\label{lem:Homeomorphism condition}
Note that a bijective map $f:(X,c) \to (Y,c')$ is a homeomorphism iff $f(c(A)) = c'(f(A)))$ for all $A \subset  X$.
\end{rem}

The following theorem and its corollaries allow us to paste together a
continuous map $f:(X,c) \to (Y,c')$ from continuous maps on sets in covers of $X$.

\begin{theorem}[\cite{Cech_1966}, Theorem 17.A.16] Let $\{U_a \mid a \in A\}$ be a locally finite cover of a closure space $(X,c)$ and let $(Y,c')$ be a closure space. If $f:X \to Y$ is a map such that the restriction of $f$ to each subspace $(c(U_a), c_{c(U_a)})$ is continuous, then $f$ is continuous.
\end{theorem}

\begin{corollary}[\cite{Cech_1966}, 17.A.18]
	\label{cor:Cont function from restriction to closed subspaces}
	Let $\{V_a\}$ be a locally finite closed cover of a closure space $(X,c)$, and let $(Y,c')$ be a closure space. If $f:X \to Y$ is a map such that the restriction of $f$ to each subspace $(V_a,c_{V_a})$ is continuous, then $f$ is continuous.
\end{corollary}

\begin{corollary}[\cite{Cech_1966}, 17.A.19]
	\label{cor:Cont function from restriction to interior cover}
	Let $\{U_a\}$ be an interior cover of a closure space $(X,c)$, and let $(Y,c')$ be a closure space. If $f:X \to Y$ is a map such that the restriction of $f$ to each subspace $(U_a,c_{U_a})$ is continuous, then $f$ is continuous. 
\end{corollary}

\subsection{Topological modifications}

It will sometimes be useful to appeal to a topological closure structure
`generated by' a given closure structure on a closure space $(X,c)$. We construct this as follows.

\begin{proposition}[\cite{Cech_1966}, 16.B.1-16.B.3]
	\label{prop:Topmod}
Let $(X,c_X)$ be a closure space. Let $\tau c_X:\mathcal{P}(X) \to \mathcal{P}(X)$ denote the map 
\begin{equation*}
\tau c_X(A) = \bigcap \{ F \subset X \mid c_X(F) = F \text{ and }A \subset F \}
\end{equation*}
Then $\tau c_X$ is the finest topological closure coarser than $c_X$, i.e. $(\tau c_X)^2 = \tau c_X$, $c_X < \tau c_X$, and for any other closure operator $c$ on $X$ with $c^2 = c$ and $c_X < c$, then $\tau c_X < c$ as well.
\end{proposition}

\begin{definition}
We call $\tau c_X$ the \emph{topological modification of $c_X$}
\end{definition}

In the next proposition, we see that topological modification of a closure structure has a characterization in terms of continuous maps.

\begin{proposition}[\cite{Cech_1966}, 16.B.4]
	\label{prop:Topological modification continuity} Let $(X,c)$ be a closure space, and let $c'$ be a closure operator on $X$. $c'$ is the topological modification $\tau c$ of $c$ 
iff the following two statements are satisfied:
\begin{enumerate}
	\item $c'$ is a topological closure operator
	\item For any topological space $(Y,\tau)$, a map $f:X \to Y$ is continuous as a map from $(X,c) \to (Y,\tau)$ iff $f$ is continuous as a map from $(X,c') \to (Y,\tau)$. 
\end{enumerate}
\end{proposition}

In categorical language, this gives the following. 

\begin{corollary}
	$\tau:\catname{Cl} \to \catname{Top}$ is a functor, and is a left-adjoint to the inclusion functor $\iota:\catname{Top} \to \catname{Cl}$.
\end{corollary}

\begin{proof}
	On objects, we define $\tau(X,c) \coloneqq (X,\tau c)$. For each
	continuous function $f:(X,c) \to (Y,c')$, let $\tau f(x) = f(x)$. We
	first claim $\tau f = f$ is continuous as a function from $(X,\tau c)
	\to (Y,\tau c')$. Let $id_Y:(Y,c') \to (Y,\tau c')$ be the identity on
	$Y$. Note that $id_Y$ is continuous by Proposition \ref{prop:Identity
	cont}, since $c' < \tau c'$ by Proposition \ref{prop:Topmod}. We now have that the
	composition $id_Y \circ f: (X,c) \to  (Y,\tau c')$ is continuous by
	Proposition \ref{prop:Composition of continuous is continuous}.
	However, by Proposition \ref{prop:Topological modification continuity},
	it follows that the map $\tau (id_Y \circ f):(X,\tau c) \to (Y,\tau
	c')$ is continuous. Since $\tau (id_Y \circ f(x)) = id_Y \circ f(x) =
	f(x)= \tau f(x)$ for all $x \in X$, it follwos that $\tau f:(X,\tau c) \to (Y,\tau c')$
	is continuous, which proves the claim. 
	
	Since $\tau$ is associative, preserves composition, and preserves identity maps, we have that $\tau:\catname{Cl} \to \catname{Top}$ is a functor.
	
	To see that $\tau$ and $\iota$ are adjoints, we must show that there
        exists a natural bijection between $\Hom_{\catname{Top}}((X,\tau
        c),(Y,\tau_Y))$ and $\Hom_{\catname{Cl}} ((X, c),(Y,\tau_Y))$. However, this follows immediately from Proposition \ref{prop:Topological modification continuity}.
\end{proof}

\subsection{Connectivity}

In the discussion of covering spaces in Subsection \ref{subsec:Covering spaces}, we will need the following definitions and results on connectivity of closure spaces.

\begin{definition}
	Let $(X,c)$ be a closure space. Two subsets $A,B \subset X$ are said to be \emph{semi-separated} if there exist neighborhoods $U$ of $A$ and $V$ of $B$ such that $U \cap B = \emptyset = A \cap V$. Two subsets $A,B \subset X$ are said to be \emph{separated} iff there exist neighborhoods $U$ of $A$ and $V$ of $B$ such that $U \cap V = \emptyset$.
\end{definition}

\begin{proposition}[\cite{Cech_1966}, 20.A.3]
	Let $(X,c)$ be a closure space. Two closed subsets are semi-separated iff they are disjoint, and two open subsets are separated iff they are disjoint.
\end{proposition}

\begin{proposition}[\cite{Cech_1966}, 20.A.6 (a)]
	Two subsets $A,B \subset X$ of a closure space $(X,c)$ are semi-separated iff $(c(A) \cap B) \cup (A \cap c(B)) = \emptyset$
\end{proposition} 

\begin{definition}
	Let $(X,c)$ be a closure space. A subset $A\subset X$ is said to be
	\emph{connected} in $X$ if $A$ is not the union of two non-empty
	semi-separated subsets of $X$, i.e. if $A$ is connected and $A = U \cup
	V$ where $(c(U) \cap V) \cup (U \cap c(V))= \emptyset$, then either $U = \emptyset$ or $V = \emptyset$. 
\end{definition}

\begin{theorem}[\cite{Cech_1966}, Theorem 22.B.2]
\label{thm:Connectedness}
	A closure space $(X,c)$ is connected iff $X$ is not the union of two disjoint non-empty open subsets. A subset $A \subset X$ is connected iff $(A,c_A) \subset (X,c)$ is connected as a closure space with the subspace closure structure $c_A$.
\end{theorem}

\subsection{New structures from old} 

We end this section with the constructions of closure structures on subspaces, disjoint unions, products, and quotients, as well as some basic results. 

 \begin{definition} 
 \begin{enumerate}[itemindent=*, leftmargin=*]
 \item Let $(X,c)$ be a closure space. We define the closure operator $c_A$ on the subspace $A\subseteq X$ by $c_A(B) = \{A \cap c(B)\}$, where $B\subseteq A \subseteq X$. 
 
 \item Let $(X_\alpha,c_\alpha)_{\alpha\in I}$ be a collection of closure spaces. The closure operator $c$ on the disjoint union $\sqcup_{i\in I} (X_i,c_i)$ is defined by $c(\sqcup_{i\in I} A_i) \coloneqq \sqcup_{i\in I} c_i(A_i)$ for any subset $\sqcup_{i\in I}A_i \subseteq \sqcup_{i\in I} X_i$.
 
 \item Let $(X_\alpha,c_\alpha)_{\alpha\in I}$ be a collection of closure spaces, let $\prod_{\alpha\in I} X_\alpha$ be the Cartesian product of the underlying sets, and let $\pi_\alpha:\prod_{\alpha\in I} X_\alpha \to (X_\alpha,c_\alpha)$ be the projection mappings from $\prod_{\alpha\in I} X_\alpha$ to $X_\alpha$. For each $x \in \prod_{\alpha\in I} X_\alpha$, let 
 \[ \mathcal{V}_x \coloneqq \{ \pi_{\alpha}^{-1}(V) \mid \alpha \in I, V \subset X_\alpha \text{ a neighborhood of } \pi_{\alpha}(x) \in X_{\alpha}\}.\] 
 By Corollary \ref{cor:Closure with prescribed base and subbase}, Part \ref{item:Closure with prescribed subbase}, there exists a unique closure structure $c_\Pi$ on $\prod_{\alpha\in I} X_\alpha$ such that $\mathcal{V}_x$ is a subbase for each $x \in \prod_{\alpha\in I} X_\alpha$. We define the \emph{product closure operator} on $\prod_{\alpha\in I} X_\alpha$ to be the closure structure $c_\Pi$.
 
 \item \label{prop:Quotient closure formula} Let $f:(X,c_X) \to Y$ be an onto
	 mapping from a closure space $(X,c_X)$ to a set $Y$. For any $V\subset
	 Y$, we define the \emph{quotient closure operator}
 		\begin{equation*}
 		c_f(V) = f(c_X(f^{-1}(V))).
 		\end{equation*} 
 \end{enumerate}
 \end{definition}

Neighborhoods in subspaces have the following useful property.

\begin{proposition}
	\label{prop:Neighborhood in neighborhood is neighborhood}
	Let $(X,c_X)$ be a closure space, and suppose that $U$ is a neighborhood
	of $x \in X$. Let $V \subset U$ be a neighborhood of $x$ in the
	subspace $(U,c_U)$, where $c_U(A) = U \cap c_X(A)$ is
	the subspace closure operator. Then $V$ is a neighborhood of $x$ in
	$(X,c_X)$
\end{proposition}
\begin{proof}
	First, note that $c_X(X - V) = c_X(X - U) \cup c_X(U-V)$. Since $U$ is
	a neighborhood of $x$, $x \notin c_X(X - U)$. Since $V$ is a
	neighborhood of $x$ in $(U,c_U)$, $x \notin U \cap
	c_X(U-V)= c_U(U-V)$. However, since $x \in U$, this implies that $x \notin
	c_X(U-V)$. Therefore $x \in X - c_X(X-V)$, and $V$ is a neighborhood of
	$x$ in $(X,c_X)$.
\end{proof}

We state the following useful results about these constructions from \cite{Cech_1966}.

\begin{theorem}[\cite{Cech_1966}, Theorem 17.A.13]
\label{thm:Restriction continuous}
Let $f:(X,c_X) \to (Y,c_Y)$ be continuous, and suppose $f(X) \subset B \subset Y$. Then $f$ is continuous as a map $f:(X,c_X) \to (B,c_B)$, where $c_B$ is the subspace closure structure on $B \subset Y$ induced from $c_Y$.
\end{theorem}

We now shift our attention to results on products.

\begin{theorem}[\cite{Cech_1966}, Theorem 17.C.3]
	\label{thm:Product local base and local subbase}
	Let $x = (x_a)_{a \in A}$ be a point in the product space $X = \Pi_{a \in A} (X_a,c_a)$ and let $\{\mathcal{U}_a\}$ be a local subbase for $x_a \in X_a$ for each $a \in A$. Then the collection of sets of the form $\pi_a^{-1}(U)$, $U \in \mathcal{U}_a$ is a local subbase for $x$ in $X$. Similarly, if each $\mathcal{U}_a$ is a local base at $x_a \in X_a$, then the collection of all sets of the form $\cap_{i = 1}^n \pi^{-1}_{a_i}(U_{a_i})$, where $U_{a_i} \in \mathcal{U}_{a_i}$ for each $i \in \{1,\dots,n\}$ is a local base at $x\in X$. 
\end{theorem}

\begin{theorem}[\cite{Cech_1966}, Theorem 17.C.6]
Let $\pi_{\alpha}:\Pi_{\alpha \in A} X_{\alpha} \to X_\alpha$ be the projections of a product space onto its coordinate spaces. Then the $\pi_{\alpha}$ are continuous for each $\alpha \in A$, where $\Pi_{\alpha \in A} X_{\alpha}$ is endowed with the product closure structure $c_\Pi$.

Moreover, $c_\Pi$ is the coarsest closure on $\Pi_{\alpha \in A} X_{\alpha}$ which makes all of the $\pi_{\alpha}$ continuous.
\end{theorem}

As in the topological case, we have

\begin{theorem}[\cite{Cech_1966}, Theorem 17.C.10]
Let $\{(X_\alpha,c_\alpha)\}_{\alpha \in A}$, be a collection of closure spaces indexed by a set $A$. A map $f:(Y,c) \to (\Pi_{\alpha \in A} X_\alpha,c_\Pi)$ is continuous iff each map
\begin{equation*}
\pi_{\alpha} \circ f:(Y,c) \to (X_\alpha,c_\alpha)
\end{equation*}
is continuous for each $\alpha \in A$.
\end{theorem}		

Finally, in the following proposition we see that taking subspaces and taking
products commute, which we will need when we develop homotopy.

\begin{proposition}[\cite{Cech_1966}, 17.c.5]
\label{prop:Product of subspaces}
If $\{(X_{\alpha},c_{\alpha})\}$, $\alpha \in A$, is a family of closure spaces and
$(U_{\alpha},c_{U_\alpha})$ is a subspace of $(X_{\alpha},c_{\alpha})$
for each $\alpha \in A$, then $(\Pi_{\alpha} U_{\alpha},c_{\Pi_U})$ is a subspace of
$(\Pi_{\alpha} X_{\alpha},c_{\Pi_X})$, where $c_{\Pi_U}$ and $c_{\Pi_X}$ are the
product closure structures on $\Pi_{\alpha} U_{\alpha}$ and $\Pi_{\alpha}
X_{\alpha}$, respectively.
\end{proposition}
	
For quotients, we will need the following

\begin{proposition}
  Let $p:(X,c_X) \to Y$ be a surjective map from a closure space $(X,c_X)$ to a set $Y$. Then the quotient closure structure $c_p$ on $Y$ induced by $p$ is the finest closure structure on $Y$ which makes $p$ continuous.	
\end{proposition}
\begin{proof}
Suppose there exists a closure structure $c$ on $Y$ which is finer than $c_p$ and for which $p:(X,c_X) \to (Y,c)$ is continuous. Then there exists a set $V \subset Y$ such that $c(V) = c(p(p^{-1}(V))) \subsetneq c_p(V) = p(c_p(p^{-1}(V)))$. However, this implies that $p$ is not continuous, a contradiction. Therefore, $c_p$ is the finest closure structure on $Y$ which makes $p$ continuous, as desired. 
\end{proof}

In the following proposition and example, we show that products in $\catname{Top}$ are also products in $\catname{Cl}$, but that quotients in $\catname{Top}$ are not necessarily quotients in $\catname{Cl}$.

\begin{proposition}[\cite{Cech_1966}, Theorem 17.C.4]
\label{prop:Product of top is top} 
The product (in $\catname{Cl}$) of a family $\mathcal{T} = \{(T_a,\tau_a) \mid a \in A\}$ of topological spaces is the product (in $\catname{Top}$) of the family $\mathcal{T}$.
\end{proposition}
\begin{proof}
The proposition follows immediately from the fact that the inclusion functor $\iota:\catname{Top} \to \catname{Cl}$ is a right-adjoint, and therefore preserves limits. An alternate proof is given in \cite{Cech_1966}.
\end{proof}

\begin{example}[\cite{Cech_1966}, Introduction to Section 33.B]
We give a simple example here to show that the quotient closure structure on
the quotient of topological spaces need not be topological. Consider the four
point space $A = (1,2,3,4)$, with the closure structure
\begin{align*}
c(1) &= \{1,2\}\\
c(2) &= \{2\}\\
c(3) &= \{3,4\}\\
c(4) &= \{4\},
\end{align*}
and note that $c^2 = c$, so $(A,c)$ is a topological space. Now consider the mapping $f:\{1,2,3,4\} \to \{x_1,x_2,x_3\}$ given by 
\begin{align*}
f(i) = \begin{cases} x_i & i = 1,2\\
x_{i-1} & i = 3,4.
\end{cases}
\end{align*}
By Proposition \ref{prop:Quotient closure formula}, the quotient closure structure $c_f$ is
\begin{align*}
c_f(x_1) &= \{x_1,x_2\}\\
c_f(x_2) &= \{x_2,x_3\}\\
c_f(x_3) &= \{x_3\}.
\end{align*}
The structure $c_f$ is not topological, however, since 
\begin{equation*}
c_f^2(x_1) = \{x_1,x_2,x_3\} \neq \{x_1,x_2\} = c_f(x_1)
\end{equation*}  
\end{example}
\section{A family of closure structures on metric spaces}
\label{sec:ScMet}
We now define the closure operators on extended pseudometric spaces which will be our main examples. 

\begin{definition}
Let $X$ be a set. A map $d:X \times X \to \R$ is called a \emph{pseudometric} iff
\begin{enumerate}
\item For all $x \in X$, $d(x,x) = 0$
\item For all $x,y \in X$, $d(x,y) = d(y,x) \geq 0$
\item For all $x,y,z \in X$, $d(x,y) \leq d(x,z) + d(z,y)$.
\end{enumerate}
If, in addition, $d(x,y) = 0 \iff x = y$, then $d$ is called a \emph{metric}.
If the range of $d$ is $[0,+\infty]$, then we call $d$ an \emph{extended (pseudo)metric}.
\end{definition}

\begin{definition}
\label{def:r-PrTop} Let $(X,d_X)$ be a metric space, $x \in X$, and $A \subset X$. Let $d(x,A) \coloneqq \inf_{y\in A} d(x,y)$. Given an $r\geq 0$, we define $c_r:\mathcal{P}(X) \to \mathcal{P}(X)$ by $c_r(A) = \{x \in X \mid d(x,A) \leq r\}$.
\end{definition}

\begin{lemma}
For any $r \geq 0$, $c_{r}$ defined above is a closure operator on the metric space $(X,d_X)$, and $c_0$ is the topological closure operator on $X$ for the topology induced by the metric.
\end{lemma}
\begin{proof} The proof follows easily from the definitions. First, we see that $c_r(\emptyset) = \emptyset$ for all $r$. Next, we note that since $r \geq 0$, $A \subseteq c(A)$. Finally, note that $d(x,A_0\cup A_1) \leq d(x,A_i)$, $i = 0,1$, for any $x \in X$. Suppose that $d(x,A_0) > d(x,A_0\cup A_1)$. Then $d(x,A_1) = d(x,A_0 \cup A_1)$, and vice-versa if $d(x,B) > d(x,A\cup B)$. Therefore $d(x,A \cup B) = \min(d(x,A),d(x,B))$, so $x \in c_r(A) \cup c_r(B) \iff x \in c_r(A\cup B)$, and we conclude that $c_r(A\cup B) = c_r(A) \cup c_r(B)$ for any $r\geq 0$. 
	
To see the last statement, note that $c_0^2 = c_0$, so $c_0$ is topological, and furthermore that $x \in c_0(A) \iff d(x,A) = 0$, so $c_0$ equals the topological closure structure induced by the metric.
\end{proof}

\begin{definition} For fixed $q,r > 0$, we say that a function between metric spaces
\[f:(X,d_X) \to (Y,d_Y)\] is \emph{$(q,r)$-continuous} if, for every $\epsilon > 0$ and $x \in X$, there exists a $\delta(x) > 0$ such that \[d_X(x,x') < q+ \delta(x) \implies d_Y(f(x),f(x')) < r+\epsilon.\] An $(q,r)$-continuous function with $\delta$ independent of $x$ is called \emph{absolutely $(q,r)$-continuous}, and if $q=r$, we simply say that $f$ is $\emph{r-continuous}$.
\end{definition}

We now show that $(q,r)$-continuity on metric spaces and continuity for maps between the associated closure spaces is equivalent.

\begin{proposition}
\label{prop:Equivalence of continuity}
Let $(X,d_X)$ and $(Y,d_Y)$ be metric spaces with closure operators $c_q$ and
$c_r$, respectively. Then a map $f:X \to Y$ is continuous as a map between
closure spaces $(X,c_q) \to (Y,c_r)$ iff it is $(q,r)$-continuous.
\end{proposition}

\begin{proof} Suppose that $f:(X,c_q) \to (Y,c_r)$ is continuous as a map
	between closure spaces. For some $\epsilon > 0$, consider the set
	$V=B_{r+\epsilon}(f(x))$, the (topologically) open ball of radius $r +
	\epsilon$ centered at $f(x)$. By construction, $f(x) \in Y - c_r(Y-V)$,
	so $V$ is a neighborhood of $f(x)$. Since $f$ is continuous, by Theorem
	\ref{thm:Neighborhood continuity} there exists a neighborhood $U$ of
	$x$ with $f(U) \subset V$. Then, by definition, $x \in X - c_q(X-U)$,
	and therefore $d(x,X-U) > q$. Taking $\delta(\epsilon,x)$ so that $q + \delta(\epsilon,x) < d(x,X-U)$, we see that $B_{q+\delta(\epsilon,x)}(x) \subset U$, and therefore that $f$ is $(q,r)$-continuous.

Conversely, let $f$ be $(q,r)$-continuous, and suppose that $x \in c_q(A)$ for
some $A \subset X$. We claim that $f(x) \subset c_q(f(A))$. To see this, first
choose an $\epsilon>0$ and let $\delta(\epsilon,x) > 0$ be such
that $d_X(x,y) < q + \delta(\epsilon,x) \implies d_Y(f(x),f(y)) < r  +
\epsilon$. Since $x \in c_q(A)$, we have that $d_X(x,A) \leq q$, and therefore
there exists a $y \in A$ with $d(x,y) \leq q + \delta(\epsilon,x)$. It follows
that $d(f(x),f(y)) < r + \epsilon$. Since $f(y) \in f(A)$ and $\epsilon$ is
arbitrary, it follows that $d(f(x),f(A)) \leq r$, and therefore $f(x) \in
c_r(f(A))$, and the claim is proved. 

The above holds for any $x \in A$, so we have $f(c_q(A)) \subseteq
c_r(f(A))$. Furthermore, since $A$ is arbitrary, $f$ is continuous, as desired.
\end{proof}

We now give a number of elementary results which aid in developing intuition and which will simplify several proofs in the following section.

The identity map obeys the following

\begin{lemma}
\label{lem:Identity q-r continuous}
For any metric space $(X,d_X)$, the identity map $\Id:(X,d_X) \to (X,d_X)$ is $(q,r)$-continuous if at least one of the following holds:
\begin{enumerate}[itemindent=*, leftmargin=*]
\item \label{item: Q leq R}$q \leq r$
\item \label{item: Diameter X} The diameter of $(X,d_X)$ is at most $r$.
\end{enumerate}
\end{lemma}
\begin{proof}
Note first that, if $q \leq r$, then $c_q < c_r$. Furthermore, if the diameter of $(X,d_X)$ is less than $r$, then $c_r(A) = X$ for any nonempty $A \subset X$, and therefore $c_q < c_r$ as well. The result now follows from Proposition \ref{prop:Identity cont}. 
\end{proof}

We have the following rule for composition.

\begin{lemma}
Let $X,Y$, and $Z$ be metric spaces. Suppose that $f:X \to Y$ be $(q,r)$-continuous, and $g:Y \to Z$ is $(s,t)$-continuous. If $r \leq s$, then $g\circ f:X \to Z$ is $(q,t)$-continuous.
\end{lemma}
\begin{proof}
Consider the diagram
\begin{equation*}
(X,c_q) \xrightarrow{f} (Y,c_r) \xrightarrow{id} (Y,c_s) \xrightarrow{g} (Z,c_t)
\end{equation*}
By Lemma \ref{lem:Identity q-r continuous}, the identity map in the diagram is continuous, and by Proposition \ref{prop:Composition of continuous is continuous}, the composition $g \circ id \circ f = g \circ f$ is continuous, which proves the lemma.
\end{proof}

Similarly, we have the following lemma.

\begin{lemma}
If $f: X \to Y$ is $(q,r)$-continuous and $p\leq q$ and $r \leq s$, then $f$ is $(p,s)$-continuous.
\end{lemma}

\begin{proof}
Consider the diagram
\begin{equation*}
(X,c_p) \xrightarrow{id} (X,c_q) \xrightarrow{f} (Y,c_r) \xrightarrow{id} (Y,c_s)
\end{equation*}
By Lemma \ref{lem:Identity q-r continuous}, the identity maps above are continuous, and by Proposition \ref{prop:Composition of continuous is continuous}, the composition is continuous, which proves the lemma.
\end{proof}

\begin{proposition}
Let $f:(X,d_X) \to (Y,d_Y)$ be a map between extended pseudometric spaces. If $f$ is Lipschitz with Lipschitz constant $K$, i.e. if, for all $x,y \in X$, $d_Y(f(x), f(y)) \leq K d_X(x,y)$, then $f:(X,c_r) \to (Y,c_{Kr})$ is continuous as a map between closure spaces.
\end{proposition}

\begin{proof}
Let $f$ be Lipschitz with Lipschitz constant $K$, i.e. 
\begin{equation*}
d_Y(f(x),f(y)) \leq Kd_X(x,y)
\end{equation*} for all $x,y \in X$. Therefore, for any $\epsilon > 0$, we have

\begin{equation*}
d_X(x,y) < r + \frac{\epsilon}{K} \implies d_Y(f(x),f(y)) < Kr + \epsilon.
\end{equation*}
Taking $\delta \coloneqq \frac{\epsilon}{K}$, the result follows from 
Proposition \ref{prop:Equivalence of continuity}.
\end{proof}

The following corollary follows immediately from the above with $K=1$, but is nonetheless worth stating separately.

\begin{corollary}
\label{cor:Contraction is continuous}
Let $f:(X,d_X) \to (Y,d_Y)$ be a map between extended pseudometric spaces. If $f$ is a contraction, i.e. $d_Y(f(x),f(y)) \leq d(x,y)$, then $f:(X,c_r) \to (Y,c_r)$ is continuous as a map between closure spaces.
\end{corollary}

\begin{proposition}
Let $d_X, d'_{X}$ be extended pseudometrics for a set $X$ and suppose that
$d_X(x,y) = d'_{X}(x,y)$ for all $x,y \in X$ for which $\min
\{d_X(x,y),d'_X(x,y)\} \leq r_0$, where $r_0 > 0$ is a constant independent of the points $x$ and $y$. Then $c_r=c'_r$ for any $0 \leq r < r_0$, where $c_r$ and $c'_r$ are the closure structures associated to $d$ and $d'$, respectively.
\end{proposition}

\begin{proof}
Let $A \subset X$, and let $0 \leq r < r_0$. Then $c_r(A) \coloneqq \{x \in X
\mid d(x,A) \leq r\}$ and $c'_r(A) = \{x \in X \mid d'(x,A) \leq r\}$ by
definition. Let $x \in c_r(A)$. Thereofre, for every sufficiently small
$\epsilon > 0$, there exists an $x_\epsilon \in A$ such that $d(x,x_\epsilon) <
r+ \epsilon$. By hypothesis, $d(x,x_\epsilon) = d'(x,x_\epsilon) < r+ \epsilon$
as well, and therefore $x \in c'_r(A)$, so $c_r(A) \subset c'_r(A)$. The same
argument with the roles of $c_r$ and $c'_r$ reversed shows that $c'_r(A)
\subset c_r(A)$, and therefore $c_r(A)=c'_r(A)$. Since the subset $A \subset X$ was arbitrary, $c_r = c'_r$.   
\end{proof}

The following corollary now follows immediately, where we let $d'(x,y) = \infty$ if $x \in X_\alpha$, $y \in X_\beta$, $\alpha \neq \beta$.

\begin{corollary}
Let $0\leq r < r_0$. Suppose $(X,d)$ may be partitioned into subsets $\{X_\alpha\}_{\alpha \in A}$ such that the $d(X_\alpha,X_\beta) \geq r_0$ for any $\alpha \neq \beta$. Then $(X,c_r) = \coprod_{\alpha \in A} (X_\alpha,c_r)$.
\end{corollary}

\begin{proposition}
Let $f:(X,c_r) \to (Y,c_s)$ be a map. If there exist points $x,x' \in X$ such that $d_X(x,x') \leq r$ and $d_Y(f(x),f(x')) > s$, then $f$ is discontinuous.
\end{proposition}

\begin{proof}
The result follows directly from Proposition \ref{prop:Equivalence of continuity} since, for sufficiently small $\epsilon>0$, there does not exist a $\delta>0$ which would make $f$ $(r,s)$-continuous at $x$ or $x'$.
\end{proof}

\begin{proposition}
Let $(X,d_X)$ be a metric space, and let $A \subset X$ be a subset with the induced metric. Then for all $r \geq 0$, the inclusion map $i:(A,c_r)\to (X,c_r)$ is a subspace inclusion in $\catname{Cl}$, i.e. $(A,c_r) = (A,c_A)$.
\end{proposition}

\begin{proof}
From the definitions or $(A,c_r)$ and $(A,c_A)$, we have $c_A(B) = c_r(B) \cap A$ for any $B \subset A$, and the result follows.
\end{proof}

\begin{proposition}
Let $(X, d_X)$ be a metric space, and define a graph $G_r = (V,E)$ by $V = X$ and $E = \{ (x,y) \mid d_X(x,y) \leq r\}$. Then each open set in $(X, c_r)$ is a union of connected components of $G_r$. In particular, if $G_r$ is connected, then $(X, c_r)$ is has no proper open sets.
\end{proposition}

\begin{proof}
First consider the case $r=0$. Since every point in $G_0$ is a connected
component, every subset of $G_0$ is a collection of connected components as well, and the result is true in this case. 

Now let $r>0$, suppose that $A\subset X$ is open and that $A$ is not the union of connected components. Then there is a point $x \in A$ and a point $y \in X \backslash A$ such $(x,y) \in E$. This implies that $d(x,y) \leq r$. However, $c_r(X\backslash A) = X\backslash A$, and therefore, for any $x \in A$, $d(x,X\backslash A) > r$, a contradiction. Therefore $A$ is a union of connected components in $G_r$, proving the result. 
\end{proof}

\begin{proposition} Suppose $(K,c_K)$ is a compact closure space. Then a continuous map $f:(K,c_K) \to (X, c_r)$ has bounded image.
\end{proposition}

\begin{proof} 
Note first that, for any $x \in X$ and any $\epsilon >0$, the ball
\begin{equation*}
B_x \coloneqq B(x,r+\epsilon) \coloneqq \{y \in X \mid d(x,y) < r + \epsilon\}
\end{equation*}
is a neighborhood of $x$. Fix an $\epsilon >0$, and consider the interior cover $\mathcal{V}\coloneqq \{B_x\}_{x \in X}$ of $X$. 
By Theorem \ref{thm:Neighborhood continuity}, 
$\{f^{-1}(B_x) \mid  \in \mathcal{V}\}$ is an interior cover of $K$.
We recall from Definition \ref{def:Compactness} that $K$ is compact iff every
interior cover of $K$ has a finite (not necessarily interior) subcover. Therefore, there exist a finite number of points $\{x_1, \dots,
x_n\}$ such that $\{f^{-1}(B_{x_i})\}_{i=1}^n$ is a cover of $K$. The
image of $f$ is therefore contained in $\cup_{i=1}^n B_{x_i}$ and it follows
that $Image(f)$ is bounded.
\end{proof}

\begin{proposition}
Let $f, g : X \to Y$ where $(X, c)$ is a closure space and $(Y,d)$ is a metric space. If
$f : (X, c) \to (Y, c_r)$ is continuous, and $d_\infty(f, g) \leq s$, then $g : (X, c) \to (Y, c_{r+2s})$ is
continuous.
\end{proposition}

\begin{proof}
Let $y \in g(c(A))$ for some $A \subset X$. We must show that $y \in c_{r+2s}(g(A))$. By definition, there exists an $x \in c(A)$ such that $g(x)=y$, and by hypothesis $d(f(x),g(x)) \leq s$. However, since $f$ is continuous, $f(c(A)) \subset c_r(f(A))$, and therefore $d(f(x),f(A)) \leq r$. It follows that, for any $x' \in A$,
\begin{align*}
d(g(x),g(x')) &\leq d(g(x),f(x)) + d(f(x),f(x')) + d(f(x'),g(x'))\\
&\leq d(f(x),f(x')) + 2s.
\end{align*}
Taking the infimum over $x' \in A$ on both sides, we obtain \begin{align*}
d(g(x),g(A)) \leq d(f(x),f(A)) + 2s \leq r + 2s
\end{align*} as desired.
\end{proof}

\section{Homotopy theory on closure spaces}
\label{sec:HomConv}

In this section, we present the construction of the fundamental groupoid and fundamental group of a closure space $(X,c)$, and then proceed to prove a general form of the Seifert-Van Kampen theorem in this setting, using interior covers of $(X,c)$. Finally, after briefly describing covering spaces in $\catname{Cl}$, we calculate the fundamental group of the circle for a family of closure structures from Section \ref{sec:ScMet}.

\subsection{The Fundamental Groupoid of a Closure Space}
\label{subsec:Homotopy groups}

\begin{definition}
\label{def:Groupoid product}
Let $I$ denote the interval $[0,1]$ endowed with the topological closure structure $\tau$, and let $(X,c),(Y,c')$ be closure spaces. We will write $(X,c) \times I$ for the product closure space $(X \times I, c_\Pi)$. Let $f,g:(X,c) \times I \to (Y,c')$ be continuous maps with $f(x,1) = g(x,0)$ for all $x \in X$. We define the map $f \star g:(X,c) \times I \to (Y,c')$ by
\begin{equation}
\label{eq:Groupoid product def}
f \star g(x,t) = \begin{cases}
f(x,2t) & t \in \left[0,\frac{1}{2}\right]\\
g(x,2t-1) & t \in \left[\frac{1}{2},1\right].
\end{cases}
\end{equation}
If $f,g:I \to (X,c)$ are continuous maps, then let $f',g':\{*\} \times I \to (X,c)$ be the maps $f'(*,t) = f(t)$ and $g'(*,t) = g(t)$. Abusing notation, we define $f \star g(t) = f' \star g'(*,t)$.
\end{definition}

\begin{proposition}
\label{prop:Groupoid product continuous}
Given continuous functions $f,g:(X,c) \times I \to (Y,c')$, the function $f\star g:(X,c) \times I \to (Y,c')$ defined in Definition \ref{eq:Groupoid product def} is continuous. Similarly, if $f,g: I \to (Y,c')$ are continuous, then $f \star g:I \to (Y,c')$ is continuous.
\end{proposition}
\begin{proof}
	Consider the cover $\{X_0,X_1\}$ of $X \times I$ by closed sets $\{X_0,X_1\}$ given by $X_0 = X \times \left[0,\frac{1}{2}\right]$ and $X_1 = X \times \left[\frac{1}{2},1\right]$. The proposition now follows from Corollary \ref{cor:Cont function from restriction to closed subspaces}. 
\end{proof}

We define homotopic maps and homotopy equivalence of spaces as in the topological setting.

\begin{definition}
We write $(X,A,c)$ for a \emph{closure space pair}, or simply \emph{a pair}, where $(X,c)$ is a closure space and $(A,c_A) \subset (X,c)$ a subspace with the subspace closure $c_A$. We will write $(X,A)$, omitting the closure structure, when there is no ambiguity. We denote the pair $(X \times I,A \times I,c_\Pi)$ by $(X,A,c) \times I$, where $c_\Pi$ is the product closure structure. A continuous map $f:(X,A,c) \to (Y,B,c')$ between pairs is a continuous map $f:(X,c) \to (Y,c')$ such that $f(A) \subseteq B$. 
\end{definition}

\begin{rem}
We recall that, by Proposition \ref{prop:Product of subspaces}, $(A \times I,
c_\Pi)$ is a subspace of $(X \times I,c_\Pi)$. 
\end{rem}

\begin{definition}
Let $f,g:(X,A,c) \to (Y,B,c')$ be continuous maps between pairs, and suppose $X' \subset X$. We say that \emph{$f$ is homotopic to $g$ rel $X'$}, denoted by $f \simeq g$ rel $X'$, iff there exists a continuous map $H:(X,A,c) \times I \to (Y,B,c')$ such that $H(x,0) = f(x)$ and $H(x,1) = g(x)$ for all $x \in X$, and $H(x,t) = f(x) = g(x)$ for all $x \in X', t \in I$. If $X' = \emptyset$, we simply say that \emph{$f$ is homotopic to $g$}.

Suppose now that for a function $f:(X,A,c) \to (Y,B,c')$, there exists a function $g:(Y,B,c') \to (X,A,c)$ such that $f\circ g \simeq id_Y$ and $g \circ f \simeq id_X$. Then we say that $f$ is a \emph{homotopy equivalence} between $(X,A,c)$ and $(Y,B,c')$, $g$ is the \emph{homotopy inverse} of $f$, and $(X,A,c)$ and $(Y,B,c')$ are \emph{homotopy equivalent}.   
\end{definition}

\begin{rem}
As in the topological case, homotopy and homotopy equivalence are equivalence relations. The proof follows from Proposition \ref{prop:Groupoid product continuous}, and is identical to that in the topological case.
\end{rem}

\begin{theorem}
Let $(X,c)$ be a closure space. Then the collection points of $X$ form the objects of a category, denoted $\Pi(X,c)$,whose morphisms are the set of homotopy classes of continuous maps $f:(I,\tau) \to (X,c)$ $\text{rel }\partial I$, denoted $[(I,\tau),(X,c)]$, and where composition of morphisms in the category is defined as follows: For any two elements $[v],[w] \in [(I,\tau),(X,c)]$ such that $v(1) = w(0)$, $[w]\circ[v] \coloneqq [v \star w]$. Furthermore, $\Pi(X)$ is a groupoid. 
\end{theorem}

\begin{rem}
Note that the product $[w]\circ [v]$ in $\Pi(X)$ is defined in the correct order for the categorical convention.
\end{rem}

\begin{proof}
The proof proceeds as for topological spaces (as in, for instance,
\cite{Spanier_1966}), with several additional verifications to confirm that the
homotopies involved are continuous for the respective closure structures. Let
$\Pi(X,c)$ be as in the statement of the theorem. 

We first show that the constant maps $c_x:I \to (X,A)$, $c_x(t) = x$, represent
the identities in $\Pi(X,c)$. Let $[f] \in \Pi(X)$ with $f(1) = x$. Consider the map $j:I \times I \to I$ given by 
\begin{equation*}
j(s,t) = \frac{t}{2} + \left(1-\frac{t}{2}\right)s
\end{equation*}
Since $j$ is continuous, the map $H = c_x\star f(j(s,t))$ is continuous as well, by Proposition \ref{prop:Groupoid product continuous}. However, $H(s,0) = c_x \star f(s)$, $H(s,1) = c_x \star f\mleft(\frac{1}{2} + \frac{s}{2}\mright) = f(s)$, and for all $t\in I$, $H(0,t) = c_x\star f\mleft(\frac{t}{2}\mright) = x = c_x \star f(0)$, and $H(1,t) = c_x \star f(1)$. Therefore $H$ is a homotopy rel $\partial I$, and $c_x$ is a left identity for $f$. The argument showing $f \star c_x \simeq f$ rel $\partial I$ is similar.

Let $[f] \in \Pi(X,c)$, and let $g:I \to (X,c)$ be given by $g(t) = f(1-t)$. Since $f(1) = g(0)$, $f\star g$ is defined, and it is continuous by Proposition \ref{prop:Groupoid product continuous}. Note, furthermore, that $g(1) = f(0)$. We define the function $H:I \times I \to (X,c)$ by
\begin{equation*}
H(s,t) = \begin{cases}
f(0) & 0 \leq s \leq \frac{t}{2}\\
f(2s-t) & \frac{t}{2} \leq s \leq \frac{1}{2}\\
f(2-2s-t) = g(2s+t-1) & \frac{1}{2} \leq s \leq 1-\frac{t}{2}\\
f(0) & 1 - \frac{t}{2} \leq s \leq 1
\end{cases}
\end{equation*} 
$H$ is well-defined by the definitions of $f$ and $g$, and $H$ is  continuous
by Corollary \ref{cor:Cont function from restriction to closed subspaces},
since $H$ is continuous on each of the closed regions defined in the right-hand
column of its definition. Furthermore, $H(s,0) = f \star g(0)$ and $H(s,1) = f
\star g(s)$ for all $s \in I$. $H$ is therefore a homotopy between $f\star g$
and the constant function $c_{f \star g(0)}$, so $f \star g \simeq 1$ rel $\partial I$. The proof for $g \star f$ is similar. Therefore, $g = f^{-1}$.

Finally, consider the continuous functions $f:(I,\tau) \to (X,c_X)$, $g:(I,\tau) \to (X,c_X)$, and $h:(I,\tau) \to (X,c_X)$ with $f(1) = g(0)$, and $g(1) = h(0)$, and form the concatenations $(f \star g) \star h$ and $f \star (g\star h)$. Let $H:I \times I \to (X,c)$ be the function
\begin{equation*}
H(s,t) = \begin{cases}
f\mleft(\frac{4t}{s+1}\mright), & t \in \left[0,\frac{s+1}{4}\right], s\in I\\
g(4t-s-1), & t \in \left[\frac{s+1}{4},\frac{s+2}{4}\right], s\in I\\
h\left(\frac{4t-2-s}{2-s}\right), & t \in \left[\frac{s+2}{4},1 \right]
\end{cases}
\end{equation*} 
$H$ is well-defined by the conditions on the endpoints of $f,g$ and $h$, and $H$ is  continuous by Corollary \ref{cor:Cont function from restriction to closed subspaces}. Since $H(0,t) = (f \star g) \star h(t)$ and $H(1,t) = f \star (g \star h)(t)$, we have $(f \star g) \star h \simeq f \star (g \star h)$.

We have now shown that $\Pi(X,c)$ is a groupoid. 
\end{proof}

\begin{definition}
	We call $\Pi(X,c)$ the \emph{fundamental groupoid of $(X,c)$}, which we
	write as $\Pi(X)$ when the structure $c$ is unambiguous. 
Given a subset $A
\subset X$, we let $\Pi(X,A)$ denote the full subgroupoid of $X$ whose objects
are the points of $A$, which we call the  \emph{fundamental groupoid of
$(X,A,c)$}. We further define $\pi_1(X,*) \coloneqq \Pi(X,*)^{op}$, which we
call the \emph{fundamental group} of $(X,c)$. We will write the product in
$\pi_1(X,*)$ as $[v][w] = [v \star w]$, and we write $\pi_1(X)$ when the
basepoint is understood.
\end{definition}

\begin{rem}
We define $\pi_1(X,*)$ as $\Pi(X,*)^{op}$ and not as $\Pi(X,*)$ in order to
make the product in the fundamental group agree with the classical definition,
i.e. where the group product of classes in $\pi_1(X)$ is written in the same order as the $\star$ product on functions, although clearly $\Pi(X,*)$ and $\Pi(X,*)^{op}$ are isomorphic groups.
\end{rem}

\subsection{A Seifert-van Kampen Theorem}

Using the above lemmas, we present a Seifert-van Kampen Theorem for the
groupoids $\Pi(X,X_0)$ using general interior covers, based on the proofs in
\cites{Brown_Salleh_1984, Brown_1982, Crowell_1959}. We begin by proving an important lemma that will be used in many of our computations.

\begin{lemma}
\label{cor:Lebesgue lemma for cubes}
Let $f:(I^n,\tau) \to (X,c_X)$ be continuous, and let $\mathcal{V}$ be an
interior cover of $X$. Then there exists a decomposition of $I^n$ into a finite number of smaller cubes such that the image of $f$ restricted to each cube is contained in a set $V \in \mathcal{V}$.
\end{lemma}

\begin{proof}
First, it follows from Theorem \ref{thm:Neighborhood continuity} that
$f^{-1}(\mathcal{V})$ is an interior cover of $I^n$. Let $\mathcal{U} = \{
\Int(f^{-1}(V)) \mid V \in \mathcal{V} \}$. Note that every set in $\mathcal{U}$
is open in $I^n$, making $\mathcal{U}$ an open cover of $I^n$. Let $\lambda$ be
the Lebesgue number for $\mathcal{U}$, and choose $k\in \mathbb{N}$ which satisfies $\frac{1}{k} < \frac{2\lambda}{\sqrt{n}}$. Decomposing $I^n$ into $k^n$ cubes whose side is length $\frac{1}{k}$ now gives the result.
\end{proof}

\begin{theorem}
\label{thm:General groupoid VK}
Let $\mathcal{U}\colonequals\{U_\alpha\}_{\alpha \in A}$ be an interior cover of a closure space $(X,c)$.  For every $\nu= (\nu_1, \dots, \nu_n) \in A^n$, let $U^{\nu} = \cap_{i=1}^n U^{\nu_i}$. Let $X_0 \subset X$, and denote $U^\nu_0 \colonequals X_0 \cap U^\nu$. Consider the diagram
\begin{equation}
\label{diag:General groupoid van Kampen}
\begin{tikzcd}
\bigsqcup_{\nu \in A^2} \Pi(U^\nu,U^\nu_0) \ar [r,shift left,"a"] \ar [r,shift right,"b",swap] & \bigsqcup_{\alpha \in A} \Pi(U^\alpha,U^\alpha_0) \ar [r,"c"] &\Pi(X,X_0) 
\end{tikzcd}
\end{equation}
where the maps $a$ and $b$ are determined by the inclusions
\begin{equation*}
a_\nu: U^{\nu_1} \cap U^{\nu_2} \to U^{\nu_1}, \text{\quad} b_\nu:U^{\nu_1} \cap U^{\nu_2} \to U^{\nu_2}.
\end{equation*}
If $X_0$ meets each path-component in the two-fold and three-fold intersections of distinct sets in $\mathcal{U}$, then $c$ is the coequalizer of $a,b$ in the category of groupoids. 
\end{theorem}

The proof of the theorem closely follows the proof given in
\cite{Brown_Salleh_1984} for topological spaces, with the addition that we use
Lemma \ref{cor:Lebesgue lemma for cubes} to guarantee the existence of the required subdivisions of $I$ and $I^2$. We refer the reader to \cite{Brown_Higgins_Sivera_2011}, Section 1.4 for a detailed discussion of the ideas in the proof in the case where $\mathcal{U} = \{U_1,U_2\}$, as well as several illustrative diagrams of the homotopies involved.

We will need the following two lemmas.

\begin{lemma}
	If a pair of topological Hausdorff spaces $(X,A)$ satisfies the homotopy extension property in $\catname{Top}$, then $(X,A,\tau)$ satisfies the homotopy extension property in $\catname{Cl}$.
\end{lemma}

\begin{proof}
	In $\catname{Top}$, the homotopy extension property is equivalent to the existence of a retract $r:X \times I \to X \times \{0\} \cup A \times I$. The existence of such a retract, however, is also sufficient to prove the homotopy extension property in $\catname{Cl}$. That is, given a closure space $(Z,c_Z)$ and continuous maps $f:A \times I \to Z$ and $g:X \times \{0\} \to Z$ such that $f(a,0) = g(a,0)$ for $a \in A$, we define a map $F:X \times I \to Z$ by 
	\begin{equation*}
	F(x,t) = \begin{cases} f \circ r(x,t) & \forall (x,t) \text{ with } r(x,t) \in A \times I\\
	g \circ r(x,t) & \forall (x,t) \text{ with } r(x,t) \in X \times \{0\}.\end{cases}	\end{equation*}
	$F$ is well-defined by the assumptions on $f$ and $g$. Since $A$ and $X$ are Hausdorff and the inclusion $A \hookrightarrow X$ is a cofibration in $\catname{Top}$, $A$ must be closed in $X$, and therefore the collection $\{X\times \{0\}, A \times I\}$ is a locally finite closed cover of $X \times I$. Since the restrictions of $F$ to each of the closed subsets $A \times I$, $X \times \{0\} \subset X \times I$ are continuous by definition, it follows from Corollary \ref{cor:Cont function from restriction to closed subspaces} that $F$ is continuous.
\end{proof}

\begin{lemma}
\label{lem:New homotopy from old} Let $\mathcal{U} = \{U^\alpha\}_{\alpha \in A}$ be an interior cover of a closure space $(X,c)$, indexed by the set $A$. Let $n\geq 1$, and suppose that $f:I^n \to (X,c)$ is a map of a topological cube such that $f$ maps the set $I^n_0$ of vertices of $I^n$ to a set $X_0 \subset X$. Suppose, furthermore, that $X_0$ meets each path component of every $k$-fold intersection of distinct sets of $\mathcal{U}$ for any $2 \leq k \leq n+1$.  

Then, $I^n$ may be subdivided into sub-cubes $\{c_\lambda\}_{\lambda \in B}$ by
planes parallel to $x_i = 0$, $i = \{0, \dots, n\}$, such that, for every
$\lambda \in B$ there exists an $\alpha(\lambda) \in A$ with $f(c_\lambda)
\subset U^{\alpha(\lambda)}$ for some $\alpha \in A$. In addition, there exists a map $g: I^n \to (X,c)$ with $f \simeq g$ rel $I^n_0$ via a homotopy $H:I^n \times I \to (X,c)$, where $g$ and $H$ satisfy
\begin{enumerate}
	\item For each $\lambda \in B$, $H(c_\lambda \times I) \subset
		U^{\alpha(\lambda)}$. (Note that $H(c_\lambda
		\times \{0\}) = f(c_\lambda) \subset U^{\alpha(\lambda)}$ by hypothesis.)
	\item For every vertex $v$ of every cube $c_\lambda$, $g(v) \in X_0$.
	\item If $e \subset c_\lambda$ is a face of $c_\lambda$ and $f(e) \subset X_0$, then $H(e \times I) = f(e) \subset X_0$.
\end{enumerate}
Furthermore, the cubical subdivision can be taken to be a refinement of any
pre-specified cubical subdivision of the cube. In particular, we can arrange
for the subdivision to refine pre-specified cubical subdivisions on the faces
$I^{n-1} \times \{0\}$ and $I^{n-1} \times \{1\}$.
\end{lemma}
\begin{proof}
	Let $\mathcal{U}' \coloneqq f^{-1}(\mathcal{U})$ and note that it is an interior cover of $I^n$ by Theorem \ref{thm:Neighborhood continuity}. Since the Lebesgue covering dimension of $I^n$ is $n$, there exists an open refinement $\mathcal{V}$ of $\mathcal{U}'$ such that every point in $I^n$ meets at most $n+1$ sets of $\mathcal{V}$. By Corollary \ref{cor:Lebesgue lemma for cubes}, there exists a decomposition of $I^n$ by planes parallel to $x_i = 0$, $i = \{0, \dots, n\}$ into sub-cubes $c_\lambda$ such that each $c_\lambda \subset V$, for some $V \in \mathcal{V}$. If $I^{n-1} \times \{0\}$ and $I^{n-1} \times \{1\}$ are already given subdivisions $d'_{\lambda'}$ and $d''_{\lambda''}$, respectively, then we may further subdivide $I^n$ so that the restriction of the cubes $c_\lambda$ to  $I^{n-1} \times \{0\}$ and $I^{n-1} \times \{1\}$ refines $d'$ and $d''$, respectively. Now let $I^{(0)} \subset I^{(1)} \subset \cdots \subset I^{(n)} = I^n$ be a cell decomposition of $I^n$ where $I^{(0)}$ is the union of the vertices of the sub-cubes $c_\lambda$, and the cells in dimension $k$ are the $k$-dimensional faces of the $c_\lambda$.  We proceed by induction on the dimension $k$ of the skeleta of $I^n$.
		
	We first consider the case $k=0$. Consider a point $v \in I^{(0)}$, and let $c_{\lambda_1},\dots,c_{\lambda_m}$ be the 
	sub-cubes that contain $v$. Let $V_i \in \mathcal{V}$, $i\in \{0,\dots,m\}$ be open 
	sets in $\mathcal{V}$ where $c_{\lambda_i} \subset V_i$ for each $i$. Since each 
	$V_i$ contains $v$, it follows that $m\leq n+1$. For each $V_i$, we choose a $U^{\alpha_i} \in \mathcal{U}$ such that $f(V_i) \subset 
	U^{\alpha_i}$. 
	
	If at least two of the $U^{\alpha_i}$ are distinct, then $X_0$ intersects each path component of $U \coloneqq \cap_{i=0}^m U^{\alpha_i}$. However, we also have that $f(v) \in U$, and therefore there exists a point $x \in X_0 \cap U$ and a path $\gamma_{vx}:I \to X$ such that $\gamma_{vx}(I) \in U, \gamma(0)=v, \gamma(1) = x$. If $v \in X_0$, then we take $x = v$ and $\gamma_{vx}$ to be the constant map. 
	
	Suppose, on the other hand, that $U^{\alpha_i}=U^{\alpha_j}$ for all
	$i,j \in \{0,\dots,m\}$. We call this set $U$. Let $a:I \to I^n$ in
	$I^n$ be a path in $I^n$ with $a(0) = v$, $a(1) = (0,\dots,0)$. If $f\circ a \subset U$,
	then $f\circ a$ is a path in $U$ from $f(v)$ to $f(0,\dots,0) \in X_0$,
	and we let $\gamma = f \circ a$. 

	If $f\circ a \not\subset U$, then let
	$U' = \cup_{U'' \in \mathcal{U}, U'' \neq U} U''$. Then $U' \neq U$ by
	definition, and $\{U, U'\}$ is an interior cover of $X$. By Corollary 
	\ref{cor:Lebesgue lemma for cubes}, there is a subdivision of $I$ into intervals 
	$[t_i, t_{i+1}]$ such that each interval is mapped to either $U$ or $U'$ by 
	$f \circ a$. Since $f\circ a(0)= f(v) \in U$, then we have $f \circ a([0,t_1])
	\subset U$. Let $[t_i,t_{i+1}]$ be the first interval with $f \circ
	a([t_i,t_{i+1}]) \subset U'$. Then $f \circ a([0,t_i])$ joins $f(v) \in
	U$ to $f \circ a(t_i) \in U \cap U''$, where $U'' 
	\in \mathcal{U}$, $U'' \neq U$. Since $X_0$ intersects every path component of every two-fold intersection 
	of distinct elements of $\mathcal{U}$, it follows that there is a path $b$ joining a point 
	$w \in X_0 \cap U'' \cap U$ to $f \circ a(t_i)$. We now let $\gamma$ be the 
	concatenation of $f\circ a([0,t_i])$ and $b$, which, by construction, joins $f(v)$ with a point in $X_0$.
	
	 Define $H^{(0)}:I^{(0)} \times I\to X$ by 
	 $H^{(0)}(v,t) = \gamma_{vx}(t)$, and let $g^{(0)}:I^{(0)} \to X$ be 
	$g^{(0)}(v) \coloneqq H^{(0)}(v,1) = \gamma_{vx}(1)$. Note that $H^{(0)}(v \times I) \subset U \subset 
	U^{\alpha_i}$ for all $i \in \{1,\dots,m\}$, that $g^{(0)}(v) \in X_0$, and 
	that $H^{(0)}(I^n_0 \times I) \subset X_0$, and, in particular, that $H^{(0)}|_{I^n_0 \times I}$ 
	is constant. This completes the initial step of the induction.
		
	Now suppose that, for any $m \in \{0,\dots,k\}$, $H^{(m)}:I^{(m)} \times I \to (X,c)$ and $g^{(m)}:I^{(m)} \to (X,c)$ satisfy
	\begin{enumerate}
		\item \label{item:Induction 1} If $f(c_\lambda) \subset
			U^{\alpha(\lambda)}$, then $H^{(m)}(c^{(m)}_\lambda
			\times I) \subset U^{\alpha(\lambda)}$.
		\item \label{item:Induction 2} For every vertex $v \in c_\lambda$, $\forall \lambda \in B$, $g^{(m)}(v) \in X_0$.
		\item \label{item:Induction 3} If $e \subset c^{(m)}_\lambda$ is a face of $c^{(m)}_\lambda$ and $f(e) \subset X_0$, then $H^{(m)}(e \times I) \subset X_0$.
	\end{enumerate}
	(Note that our initial step proved these conditions for $H^{(0)}$ and $g^{(0)}$.) For each $(k+1)$-dimensional face $\mu$ in $I^{(k+1)}$, let $\mu^{(k)} \coloneqq \mu \cap I^{(k)}$, i.e. the $k$-dimensional skeleton of $\mu$. Since each $\mu \times I$ is a $(k+2)$-dimensional cube, there exists a retract $r_{\mu}:\mu \times I \to (\mu \times {0}) \cup \mu^{(k)} \times I \subset \mu \times I$, and we define $H_{\mu}^{(k+1)}:\mu \times I \to (X,c)$ by 
	\begin{equation*}
	H_{\mu}^{(k+1)}(x,t) = \begin{cases}
	H^{(k)}(r_{\mu}(x,t)), & r_{\mu}(x) \in \mu^{(k)} \times I \subset \mu \times I \\
	f(r_{\mu}(x,t)), & r_{\mu}(x,t) \in \mu \times \{0\} \subset \mu \times I.
	\end{cases}
	\end{equation*}
	First, since $H^{(k)}(x,0) = f(x) = f(r_{\mu}(x,0))$, it follows that
	$H^{(k+1)}_{\mu}$ is well-defined, and it is continuous by Corollary
	\ref{cor:Cont function from restriction to closed subspaces}.
	Furthermore, if two cubes $\mu$ and $\mu'$ share a face $\mu''$, then
	it follows from the definition that $H^{(k+1)}_\mu|_{\mu''} =
	H^{(k+1)}_{\mu'}|_{\mu''}$. We may therefore define the function
	$H^{(k+1)}:I^{(k+1)} \times I \to (X,c)$ by $H^{(k+1)}(x,t) \coloneqq
	H^{(k+1)}_{\mu}(x,t)$ for $(x,t)\in \mu$. $H^{(k+1)}$ is continuous by
	Corollary \ref{cor:Cont function from restriction to closed subspaces},
	and by construction, $H^{(k+1)}$ extends $H^{(k)}$. We let
	$g^{(k+1)}(x) \coloneqq H^{(k+1)}(x,1)$ for all $x \in X^{(k+1)}$. It
	follows directly from the induction hypothesis and the construction of
	$H^{(k+1)}$ that the functions $H^{(k+1)}$ and $g^{(k+1)}$ satisfy
	Items \ref{item:Induction 1} and \ref{item:Induction 2} above, with $H^{(k)}$ and $g^{(k)}$ replaced by $H^{(k+1)}$ and $g^{(k+1)}$, respectively.
	
	To see Item \ref{item:Induction 3}, let $\mu'\subset \mu$ be a face of some $c^{(k+1)}_\lambda$ with $f(\mu') \subset X_0$. Recall that $r_{\mu'}(x,t) \in (\mu' \times {0}) \cup (\mu')^{(k-1)} \times I$ for any $(x,t) \in \mu'$. By the induction hypothesis, $H^{(k)}((\mu')^{(k-1)} \times I) \subset X_0$, and by the construction of $H^{(k)}$ we have that $H^{(k)}(\mu' \times {0}) = f(\mu') \subset X_0$. Therefore $H^{(k)}(x,t) \in X_0$ for any $(x,t) \in \mu' \times I$. However, since $r_\mu$ is the identity on $\mu' \times I$, $H^{(k+1)}(\mu') = H^{(k)}(\mu') \subset X_0$, and Item \ref{item:Induction 3} is satisfied.
	
	To complete the proof of the lemma, we define $H(x,t) \coloneqq H^{(n)}(x,t)$ and $g(x) \coloneqq H(x,1)$.
\end{proof}

We now proceed with the proof of Theorem \ref{thm:General groupoid VK}.

\begin{proof}[Proof of Theorem \ref{thm:General groupoid VK}]We begin with a remark on notation. Throughout the proof, the product $\star$ will be reserved for the $\star$-product of two functions, $*$ will indicate the composition of two functors between groupoids, and we will write $\circ$ for the composition of two morphisms of a groupoid.

We need to show that, for any groupoid $\Gamma$ and any functor \[F:\bigsqcup_{\alpha \in A} \Pi(U^\alpha,U^\alpha_0)\to \Gamma\] such that $F*a = F*b$, then there exists a unique functor $G:\Pi(X,X_0) \to \Gamma$ such that $G * c = F$. We first remark that for any closure space $(U,c_U)$ and $U_0 \subseteq U$, a path \[\lambda:([a,b],\{a,b\}) \to (U,U_0)\] may be taken to be a representative of the class $[\lambda(tb + (1-t)a)] \in \Pi(U,U_0)$ by reparametrizing the domain. Abusing notation slightly, we will use the parametrizations interchangeably in what follows, and, in particular, we will refer to homotopy classes of the form $[\lambda(tb + (1-t)a)]$ simply as $[\lambda([a,b])]$.

 Let $\lambda:(I,\partial I) \to (X,X_0)$ be a path in $X$ with endpoints in $X_0$, and suppose that there exists a decomposition of the $I$ into subintervals $[t_{i-1},t_i]$, $i \in \{1,\dots,k\}$, where for every $i$,
 \begin{enumerate}
 \item \label{item:Item 1 in VK proof} $\lambda(t_i) \in X_0$,
 \item There exists a $U^{\alpha_i} \in \mathcal{U}$ with $\lambda([t_{i-1},t_{i}]) \subset U^{\alpha_i}$, and
 \item \label{item:Item 3 in VK proof}$U^{\alpha_{i-1}} \neq U^{\alpha_{i}}$.
 \end{enumerate}
 If a path $\lambda$ admits such a decomposition, then we say that $\lambda$ is \emph{admissible}. Denote by $\lambda_i:[t_{i-1},t_i] \to X$ the restriction of $\lambda$ to the interval $[t_{i-1},t_i]$. Denote by $[\lambda_i] \in \Pi(X,X_0)$ the class of $\lambda_i$ in the fundamental groupoid of $(X,X_0)$. Since $\lambda(t_i) \in X_0$ for all $i$, $\lambda_i$ also represents a class in $\Pi(U^{\alpha_i},U^{\alpha_i}_0)$, which we denote by $[\lambda_i]^{\alpha_i} \in \Pi(U^{\alpha_i},U^{\alpha_i}_0)$. We often write $F[\lambda]^{\alpha_i}$, $G[\lambda]^{\alpha_i}$ for $F([\lambda]^{\alpha_i})$, $G([\lambda]^{\alpha_i})$, and other functions from $\Pi(U^{\alpha_i}, U^{\alpha_i}_0)$, when there is no ambiguity. Note that $c[\lambda_i]^{\alpha_i} = [\lambda_i]$, and that, by reparametrizing, we have $[\lambda_k] \circ \cdots \circ [\lambda_1] = [\lambda] \in \Pi(X,X_0)$. Therefore, if $G$ exists, then for any admissible path $\lambda$, $G$ must satisfy 
 \begin{equation}
 \label{eq:Functor G}
 G[\lambda] = G * c[\lambda_k]^{\alpha_k}\circ \dots \circ G * c[\lambda_1]^{\alpha_1} = F[\lambda_k]^{\alpha_k} \circ \cdots \circ F[\lambda_1]^{\alpha_1}.
 \end{equation}
 The relation $F * a = F * b$ ensures that, if the image of $\lambda$ is contained in more than one set $U^\alpha$ in the cover $\mathcal{U}$, the value of $G[\lambda]$ in Equation \ref{eq:Functor G} above is independent of the choices of the $\alpha_i$. 
 
 For a general path $\alpha:(I,\partial I) \to (X,X_0)$, Lemma \ref{lem:New homotopy from old} with $n=1$ gives a homotopy from $\alpha$ to an admissible path $\lambda$. Therefore, if $F$ exists, then we must additionally have $G[\alpha] = G[\lambda]$. Since $G[\lambda]$ must be given by Equation \ref{eq:Functor G}, we conclude that, if $G$ is well-defined, it must be unique.
 
 To see that $G$ is well-defined, it remains to show that, for two admissible homotopic paths $\lambda, \lambda':(I,\partial I) \to (X,X_0)$  which are in the same homotopy class rel $\partial I$, we have 
\begin{equation*}
G[\lambda] = F[\lambda_k]^{\alpha_k}\circ \cdots \circ F[\lambda_1]^{\alpha_1} = F[\lambda'_{k'}]^{\alpha'_{k'}} \circ \cdots \circ F[\lambda'_{1}]^{\alpha'_{1}} = G[\lambda'],
\end{equation*}
where $\{[t_{i-1},t_{i}]\}_{i=0}^k$ and $\{[t'_{i-1},t'_{i}]\}_{i=0}^{k'}$ are the subdivisions of $I$ for $\lambda$ and $\lambda'$, respectively.
Let $\tilde{H}:I^2 \to (X,c)$ be a homotopy rel $\partial I$ from $\lambda$ to $\lambda'$. Since $X_0$ intersects every path component of every two-fold and three-fold intersection of distinct sets of $\mathcal{U}$, we may apply Lemma \ref{lem:New homotopy from old} with $n=2$ to obtain maps $H:I^2 \times I \to (X,X_0)$ and $g(s,t) \coloneqq H(s,t,1)$, and a subdivision of $I^2$ into squares $\{c_{ij}\}$ such that the conclusions of Lemma \ref{lem:New homotopy from old} are satisfied. In particular, the restriction of $\{c_{ij}\}$ to $I \times \{0\}$ and $I \times \{1\}$ refines the subdivisions of $I$ given by $\{[t_{i-1},t_{i}]\}_{i=0}^k$ and $\{[t'_{i-1},t'_{i}]\}_{i'=0}^{k'}$, respectively.
 
 We first show that for the paths $H(\cdot,0,1)$ and $\lambda$, we have $G[H(\cdot,0,1)] = G[\lambda]$. Let $\{0,t_1, \dots, t_{k-1},1\}$ and $\{0,t'_1, \dots, t'_{k'-1},1\}$ be the endpoints of the subdivisions of $[0,1]$ for $\lambda$ and $\lambda'$, respectively, and let $\{0,s_1,\dots, s_{l-1},1\}$ and $\{0,s'_1,\dots,s'_{l'-1},1\}$ be the endpoints of the restriction of the subdivision $c\coloneqq \{c_{ij}\}$ to $I \times \{0\}$ and $I \times \{1\}$, respectively. Since the subdivision of $I$ given by $\{0,s_1,\dots, s_{l-1},1\}$ refines the one generated by $\{0,t_1, \dots, t_{k-1},1\}$, it follows that $\{0,t_1, \dots, t_{k-1},1\} \subset \{0,s_1,\dots, s_{l-1},1\}$. By construction, $\lambda([t_{i-1},t_i]) \subset U^{\alpha_i} \in \mathcal{U}$, and therefore we have $\lambda([s_{m-1},s_{m}]) \subset U^{\alpha_i}$ for all $t_{i-1} \leq s_{m-1}, s_m \leq t_i$. By the conclusions of Lemma \ref{lem:New homotopy from old}, $H([s_{m-1},s_m], 0, I) \subset U^{\alpha_i}$ for $t_{i-1} \leq s_{m-1}, s_m \leq t_i$ as well. Therefore, $H([t_{i-1},t_i], 0, I) \subset U^{\alpha_i}$, and this implies that the concatenation of the paths $H([t_{i-1},t_i], 0, 0) = F[\lambda_i]$ and $H(t_i, 0, I)$, is homotopic to the concatenation of the paths $H(t_{i-1}, 0, I)$ and $H([t_{i-1},t_i], 0, 1)$, i.e. $\lambda_i \star H(t_i, 0, I) \simeq H(t_{i-1}, 0, I) \star H([t_{i-1},t_i], 0, 1)$. Therefore, both of these concatenations represent the same element of $\Pi(U^{\alpha_i},U^{\alpha_i}_0)$, so 
 \begin{align*}
F[H(t_i, 0, I)]^{\alpha_i} \circ F[\lambda_k]^{\alpha_i} =  F[H([t_{i-1},t_i], 0, 1)]^{\alpha_i} \circ F[H(t_{i-1}, 0, I)]^{\alpha_i}.
 \end{align*}
 Note, as well, that $H(1,0,I)$ and $H(0,0,I)$ are constant, so $[\lambda] = [\lambda]\star[H(1,0,I)]$ and $[H(I,0,1)] = [H(1,0,I)] \star [H(I,0,1)]$. Using a succession of these homotopies, each inside the relevant $U^{\alpha_i}$, we arrive at a homotopy between $\lambda \star H(1,0,I)$ and $g(\cdot, 0)\star H(I,0,1)$, and we also see that
 \begin{align*}
 G[\lambda] &= F[\lambda_k]^{\alpha_k}\circ \cdots \circ F[\lambda_1]^{\alpha_1} \\
 &= F[H([t_{k-1},1],0,1)]^{\alpha_k}\circ F[H([0,t_1],0,1)]^{\alpha_k}\\
 &= G[g(\cdot,0)].
 \end{align*} 
 
 The proof that $G[\lambda'] = G[g(\cdot,1)]$ is identical.
 
 We must now show that $G[g(\cdot,0)] = G[g(\cdot,1)]$. The proof proceeds largely as above. We first note that, since $\tilde{H}$ is a homotopy rel $\partial I$, $\lambda(0) = \tilde{H}(0,\cdot) = \lambda'(0) \in X_0$ and $\lambda(1) = \tilde{H}(1,\cdot) = \lambda'(1) \in X_0$, and in particular both are constant functions. By the construction of $H$ in Lemma \ref{lem:New homotopy from old}, the functions $H(\epsilon,\epsilon',t)$ are constant as well, where $\epsilon, \epsilon' \in \{0,1\}$ Furthermore, $H(s,0,t) = H(r(s,0,t))$, where $r$ is the retract of the face $ 0 \times I \times I$ onto $0 \times I \times 0 \cup 0 \times \partial I \times I$. Since $H$ is constant on the latter set, it follows that $H$ is constant on $0 \times I \times I$. Similarly, $H$ is constant on $1 \times I \times I$, and therefore $g(s,t) = H(s,t,1)$ is a homotopy rel $\partial I$ between $H(\cdot,0,1)$ and $H(\cdot,1,1)$.
 
Now recall that, by Lemma \ref{lem:New homotopy from old}, $g(v_1,v_2) \in X_0$ for any vertex $v = (v_1,v_2)$ of a cube $c_{ij}$. Therefore, for all $i$ and $j$, the paths $g([t_{i-1},t_{i}],t_j)$, $g(t_{i+1},[t_j,t_{j+1}])$, $g(t_{i},[t_j,t_{j+1}])$, $g([t_i,t_{i+1}],t_{j+1}) \in U^{\alpha_{ij}}$ for some $U^{\alpha_{ij}}$ with $g(c_{ij}) \in U^{\alpha_{ij}}$. Furthermore, the contatenations 
\begin{equation*}
g([t_{i-1},t_{i}],t_j) \star g(t_{i+1},[t_j,t_{j+1}]) \simeq g(t_{i},[t_j,t_{j+1}]) \star g([t_i,t_{i+1}],t_{j+1})
\end{equation*}
are homotopic in $U^{\alpha_{ij}}$. Denote by $H_{ij}$ a homotopy that connects them. Note that that the right side of Equation \ref{eq:Functor G} is constant for two paths which are homotopic through such a homotopy $H_{ij}$. Since $g(0,t)$ and $g(1,t)$ are constant, by concatenating a sequence of the homotopies $H_{ij}$, we arrive at a homotopy between $g(0,t)\star g(s,1)$ and $g(s,0) \star g(1,t)$ which keeps Equation \ref{eq:Functor G} constant. We therefore see that $G$ depends only on the homotopy class of $\lambda$. 

Since $G$ is a functor by construction, this completes the proof.
\end{proof}

As in the topological case, the traditional version of the Seifert-van Kampen theorem for groups is now a formal consequence of Theorem \ref{thm:General groupoid VK}.

\begin{theorem}[Seifert-van Kampen Theorem for Groups]
\label{thm:Group VK}
Let $(X,c_X)$ be a closure space with interior cover $\mathcal{U}\colonequals\{U_1,U_2\}$ be an interior cover. Let $U_{12} \colonequals U_1 \cap U_2$, and suppose that $U_1,U_2$, and $U_{12}$ are path-connected and the point $* \in U_{12}$. Let $i_\alpha:U_{12} \to U_\alpha$ and $j_\alpha:U_\alpha \to X$ denote the respective inclusions. Then
\begin{equation*}
\begin{tikzcd}
\pi_1(U_{12},*) \ar[d,"\pi_1(i_2)"] \ar[r,"\pi_1(i_1)"] & \pi_1(U_1,*) \ar [d,"\pi_1(j_1)"] \\
\pi_1(U_2,*) \ar [r,"\pi_1(j_2)"] & \pi_1(X,*)
\end{tikzcd}
\end{equation*}
is a pushout in the category of groups. 
\end{theorem}

\begin{proof}
Apply Theorem \ref{thm:General groupoid VK} to the cover $\mathcal{U}$ with $X_0 = {*}$.
\end{proof}

\subsection{Covering spaces and the fundamental group of $(S^1,c_r)$}
\label{subsec:Covering spaces}
In this section, we will compute the fundamental group of the circle endowed with different closure operators, as well as that of several graphs which may be seen as subspaces of $S^1$ with the induced closure operator. The main complication we encounter is that, depending on the closure structure, we are no longer guaranteed the existence of a lift to a contractible covering space $p:E \to S^1$ for maps or homotopies of maps from $S^1$ to $S^1$, and, when a lifting naively appears to exist, such a lifting starting at a given point $x_0 \in E$ may no longer necessarily be unique. We illustrate these problems with the following examples.

\begin{definition}
	Let $R \subset \R \times \R$ be the relation $(x,y) \in R$ iff $x = y
	\mod 1$, and let $S^1$ denote the resulting quotient metric space with
	quotient metric $d_{S^1}$ (i.e. $(S^1,d_{S^1})$) is the circle with
	circumference $1$. In what follows, we will denote elements of $S^1$
	by their equivalence classes, i.e. $[x] \in S^1$, where $x \in \R$.
\end{definition}

\begin{rem}
	Note that, for any $r \geq \frac{1}{2}$, the space $(S^1,c_{r})$ is indiscrete. 
\end{rem}
\begin{lemma} 
	For any $r \geq 0$, the map $p:(\R,c_r) \to (S^1,c_r)$ given by $p(x) = x \mod 1$ is continuous.
\end{lemma}

\begin{proof}
	Since $p$ is a contraction of metric spaces, it is continuous as a map
	$(\R,c_r) \to (S^1,c_r)$ by Corollary \ref{cor:Contraction is continuous}.
\end{proof}

The following two propositions demonstrate that there are genuinely new phenomena that
must be accounted for in order to develop a theory of covering maps for closure
spaces. In Proposition \ref{prop:Example 1}, we find lifts to what one expects to be a covering space
that are not unique, and in the \ref{prop:Example 2}, we construct two
homotopic maps which do have unique lifts to a candidate covering space, but
where there exists a homotopy between them on the base space which does not have a lift.

\begin{proposition}
\label{prop:Example 1}
Let $\gamma:([0,1],c_\tau) \to (S^1,c_{1})$ be the constant map $\gamma(x) =
0$, and let $p:(\R,c_{1}) \to (S^1,c_{1})$ be the map $p(x) = [x] \in S^1$.
Then $\gamma$ does not have a unique lift $F:([0,1],c_\tau) \to (\R,c_{1})$
with $F(0) = 0$.
\end{proposition}

\begin{proof}
	We will construct two non-equal lifts of $\gamma$ which take the same value
	at $0$. Define the functions $f,g: [0,1] \to \R$ by \[f(x) = 0,\] and
\begin{equation*}
g(x) = \begin{cases}
    0 & x \in \left[0,\frac{1}{2}\right],\\
    1 & x \in \left(\frac{1}{2},1\right].
\end{cases}
\end{equation*}
Since $(S^1,c_{1/2})$ is indiscrete, both $f$ and $g$ are continuous with respect to the closure structures, and that they both lift $\gamma$ to $\R$, i.e. $pf = \gamma = pg$, starting at the basepoint $f(0)=g(0)=0$, so the lift of $\gamma$ given an initial point in this case is not unique.
\end{proof}

For the second example, we will require the following lemma. 

\begin{lemma}
\label{lem:Lifts constant integer difference}
Consider the closure spaces $(\R,c_{1/2})$ and $(S^1,c_{1/2})$, and let
 $p:(\R,c_{1/2}) \to (S^1,c_{1/2})$ be the map $p(z) \to [z] \in S^1$. Let $X$
 be a connected topological space, and let $f : X \to (S^1,c_{1/2})$ be a continuous
 map. 

 Suppose that there exists a fixed $0 < \alpha < \tfrac{1}{2}$ such that, for every $x \in X$,
there is an open neighborhood
$U\subset X$ of $x$ with $f(U) \subset I_{f(x),\alpha} \coloneqq
p((z-\alpha,z+\alpha))$, where $z \in p^{-1}(f(x))$. (Note that the set
$I_{f(x),\alpha}$ is independent
of the choice of $z \in p^{-1}(f(x))$.) Then, for any two continuous lifts $F,G:
X \to (\R,c_{1/2})$ of
$f$, there exists an integer $n_{F,G} \in \Z$ such that $F(x)-G(x) = n_{F,G}$
for all $x \in X$.
\end{lemma}

\begin{proof}
     Let
     $F,G:X \to (\R,c_{1/2})$ be continuous lifts of $f$ to $(\R,c_{1/2})$,
     i.e. $p\circ F = p \circ G = f$.
    By the connectedness of $X$, it suffices to show that the integer-valued
    function $F-G$ is locally constant.
Indeed, let $x\in X$ and write $F(x) = t$ and $G(x) = t+n$.
By hypothesis, there exists a neighborhood $U \subset X$ of $x$ such that $f(U)
\subset I_{f(x),\alpha}$ for some $0< \alpha< \tfrac{1}{2}$. By continuity, for any $\epsilon>0$, the sets
\begin{align*}
    A_\epsilon &\coloneqq
    F^{-1}\left(t-\tfrac{1}{2}-\epsilon,t+\tfrac{1}{2}+\epsilon\right)\\
    B_\epsilon &\coloneqq
    G^{-1}\left(t+n-\tfrac{1}{2}-\epsilon,t+n+\tfrac{1}{2}+\epsilon\right)
\end{align*}
are neighborhoods of $x \in X$, and it follows that the set $U_{\epsilon}' \coloneqq U
\cap A_\epsilon
\cap B_\epsilon$ is also a neighborhood of $x$. 
Choosing $\epsilon>0$
sufficiently small so that
\begin{align*}
    p^{-1}\left(f(U_{\epsilon}')\right)& \cap
    \left(t-\tfrac{1}{2}-\epsilon,t+\tfrac{1}{2}+\epsilon\right) =
    \left(t-\alpha,t+\alpha\right), \text
    { and }\\
    p^{-1}\left(f(U_{\epsilon}')\right) &\cap
    \left(t + n -\tfrac{1}{2}-\epsilon,t+n+\tfrac{1}{2}+\epsilon\right) =
    \left(t+n-\alpha,t+n+\alpha\right),
\end{align*}
it follows that
\begin{align*}
    F(U_{\epsilon}') &\subset (t-\alpha,t+\alpha), \text
    { and }\\
    G(U_{\epsilon}') &\subset (t+n-\alpha,t+n+\alpha).
\end{align*}
since $F$ and $G$ are lifts of $f$. We now observe that
\begin{align*}
    n = F(x)-G(x)&\subset (F-G)(U_{\epsilon}') = \{F(x) - G(x)\mid x \in U_{\epsilon}'\}\\ 
                                       &\subset \Z \cap \{ y - y' \mid y \in
                                       F(U_{\epsilon}'), y' \in
                                   G(U_{\epsilon}')\}\\
                                       &= \{n\},
\end{align*}
so $F-G$ must take the constant value $n$ on $U_{\epsilon}'$. Since the point $x \in X$ was
arbitrary, $U_{\epsilon}'$ is a neighborhood of $x$, and $X$ is connected, we may take $n_{F,G} = n$, and the result follows.
\end{proof}

\begin{proposition}
\label{prop:Example 2}
 Consider the closure spaces $(\R,c_{1/2})$ and $(S^1,c_{1/2})$, and let
 $p:(\R,c_{1/2}) \to (S^1,c_{1/2})$ be the map $p(x) \to [x] \in S^1$. Then
 there exist homotopic maps $f,g:([0,1],c_\tau) \to (S^1,c_{1/2})$ with unique
 lifts $\tilde{f},\tilde{g}:([0,1],c_\tau) \to (\R,c_{1/2})$ such that
 $\tilde{f}(0) = \tilde{g}(0) = 0 \in \R$, and a homotopy
 $H:([0,1]^2,c_\tau) \to (S^1,c_{1/2})$ between $f$ and $g$ such that $H$ does not have a continuous lift to a homotopy between
 $\tilde{f}$ and $\tilde{g}$
 \end{proposition}
 \begin{proof}
	 We will construct $f,g,\tilde{f},\tilde{g}$, and $H$, and show that they have the desired
properties. First, consider the maps $f,g: ([0,1],c_\tau) \to (S^1,c_{1/2})$ given by
\begin{align*}
f(x) &= [0]\\
g(x) &= \begin{cases}
    [0] & x \in \left[\tfrac{5}{6},1\right] \cup \left[0,\tfrac{1}{6}\right) \\
\left[\tfrac{1}{3}\right] & x \in \left[\tfrac{1}{6},\tfrac{1}{2}\right)\\
\left[\tfrac{2}{3}\right] & x \in \left[\tfrac{1}{2},\tfrac{5}{6}\right).
\end{cases}
\end{align*}
We first show that these two maps are homotopic given the closure structure
$c_{1/2}$. Let $H:(S^1,\tau) \times ([0,1],_\tau) \to (S^1,c_{1/2})$ be the function
\begin{align*}
H(x,t) = \begin{cases}
[0] & t \in \left[0,\frac{1}{2}\right] \\
g(x) & t \in \left(\frac{1}{2},1\right].\\
\end{cases}
\end{align*}
Since $(S^1,c_{1/2})$ is indiscrete, it follows that $H$ is continuous, and
therefore $H$ is a homotopy between $f$ and $g$. 

Now let $\tilde{f}:([0,1],c_\tau) \to (\R,c_{1/2})$ be the constant map
$\tilde{f}(x)=0, x\in [0,1]$, and let $\tilde{g}:([0,1],c_\tau) \to (\R,c_{1/2})$ be given by
\begin{align*}
\tilde{g}(x)
&= \begin{cases}
0 & x \in \left[0,\frac{1}{6}\right)\\
\frac{1}{3} & x \in \left[\frac{1}{6},\frac{1}{2}\right) \\
\frac{2}{3} & x \in \left[\frac{1}{2},\frac{5}{6}\right) \\
1 & x \in \left[\frac{5}{6},1\right],
\end{cases}
\end{align*}
Note that $\tilde{f}$ and $\tilde{g}$ are continuous lifts of $f$ and $g$,
respectively, and, since $f$ and $g$ both satisfy the hypothesis of Lemma
\ref{lem:Lifts constant integer difference} with $\tfrac{1}{3} < \alpha <
\tfrac{1}{2}$, it follows that $\tilde{f}$ and
$\tilde{g}$ are
the unique lifts of $f$ and $g$ with $\tilde{f}(0) = \tilde{g}(0) = 0$. 

Now suppose that there exists a continuous lift $\tilde{H}$ of $H$ to
$(\R,c_{1/2})$ which is a homotopy between $\tilde{f}$ and $\tilde{g}$. We first remark that, for any closed
path $\gamma:[0,1] \to [0,1]^2$, i.e. a map $\gamma$ such that $\gamma(0) = \gamma(1)$, the lift
$\tilde{H} \circ \gamma$ of
$H\circ \gamma$ must be a closed path in $(\R,c_{1/2})$, i.e.
$\tilde{H}\circ \gamma(0) = \tilde{H} \circ \gamma(1)$. Consider the path
$\gamma:[0,1] \to [0,1]^2$ defined by 
\begin{equation*}
\gamma(t) = \begin{cases} (0,4t) & t\in \left[0,\frac{1}{4}\right)\\
        \left(4(t-1/4),1\right)  & t \in \left[\frac{1}{4},\frac{1}{2}\right) \\
                \left(1,1-4(t-1/2)\right) & t \in \left[\tfrac{1}{2},\frac{3}{4}\right) \\
                \left(1-4(t-3/4),0\right) & t\in \left[\frac{3}{4},1\right].
    \end{cases}
\end{equation*}
Note that $\gamma$ is continuous by definition, being the clockwise loop around the
boundary of $[0,1]^2$, and $H\circ \gamma$ takes on
the values $[0],[1/3],[2/3],[0]$ as $t$ increases from $0$ to
$1$. Define the
map $G:[0,1] \to (\R,c_{1/2})$ by
\begin{equation*}
    G(t) \coloneqq \begin{cases} 0 & t \in \left[0,\frac{1}{4}\right) \\
        \tilde{g}(4(t-1/4)) & t \in \left[\frac{1}{4},\frac{1}{2}\right)\\
        1 & t \in \left[\frac{1}{2},1\right],
    \end{cases}
\end{equation*}
and note that $G$ is continuous by Corollary \ref{cor:Cont function from
restriction to closed subspaces}. Furthermore, we see that $G$ lifts $H\circ
\gamma$, and
since $H\circ \gamma$ satisfies the hypothesis of Lemma \ref{lem:Lifts constant
integer difference} with $\tfrac{1}{3} < \alpha < \tfrac{1}{2}$, $G$ is the unique lift of $H\circ \gamma$ with $G(0) = 0$.
However, since $\tilde{H}\circ \gamma$ is also a lift of $H \circ \gamma$ with
$0 = \tilde{f}(0) = \tilde{H} \circ \gamma(0)$, it follwos that $G =
\tilde{H}\circ \gamma$. We therefore have $\tilde{H}\circ
\gamma(0) = G(0) \neq G(1)= \tilde{H}\circ
\gamma(1)$, a contradiction. It follows that $H$ has no continuous lift to a
homotopy between $\tilde{f}$ and $\tilde{g}$.
\end{proof}

The above examples illustrate several ways in which the standard methods for computing the fundamental group of $S^1$ in the topological category do not immediately generalize to closure structures on $S^1$. We will solve these problems by incorporating the closure structure, via neighborhood systems,  explicitly into the definitions of a covering map, analogously to how the topology of a space appears in the standard definition of a covering map for topological spaces. Once the new definitions are in place, we will see that the familiar construction of the fundamental group of $S^1$ then generalizes to this setting.

We begin with a brief discussion of covering maps in closure spaces.

\begin{definition}
Let $p:(E,c_E) \to (B,c_B)$ be a surjective continuous function between closure spaces, and let $U \subset B$ be a neighborhood of a point $x \in U$. If $(F,c_F)$ is another closure space, let $q_1: U \times F \to U$ denote the projection onto the first factor. We say that a homeomorphism $\phi:p^{-1}(U) \to U \times F$ such that $q_1 \circ \phi = p$ is a \emph{trivialization of $p$ over $U$}. When such a $\phi$ exists, $p$ is said to be \emph{trivial over U}. We say that $p$ is \emph{locally trivial} if there exists an interior cover $\mathcal{U}$ of $B$ such that $p$ is trivial over $U$ for each $U \in \mathcal{U}$. 
\end{definition}

\begin{definition}
A locally trivial map $p:(E,c_E) \to (B,c_B)$ is said to be a \emph{covering
map} if, $\forall b \in B$, $F_b = p^{-1}(b)$ is discrete in the subspace closure structure on $F_b \subset E$ (i.e. every point is both open and closed), and there is an interior cover $\mathcal{U}_B$ for $B$ such that $p$ is trivial over $U$ for every $U\in \mathcal{U}_B$. We call $F_b$ the \emph{fiber of $p$ at $b$}.
\end{definition}

\begin{rem}
We observe that, for a covering $p$, since $p^{-1}(b)$ is discrete, $U 
\times p^{-1}(b)$ is homeomorphic to a disjoint union of copies of $U$ embedded in $E$. The summands $U \times \{x\}$ of $p^{-1}(U) = \sqcup_{x \in p^{-1}(b)} U \times \{x\}$ are called the \emph{sheets} of the covering $p$ over $U$.
\end{rem}

We now give a proposition, which gives sufficient conditions for the uniqueness
of a lift.

\begin{proposition}
	\label{prop:Unique liftings intermediate}
	Let $(E,c_E,*), (B,c_B,*)$, and $(X,c_X,*)$ be closure spaces with
	basepoint, suppose $(X,c_X,*)$ is
	connected, and suppose that the following diagram of continuous
	basepoint-preserving maps commutes, where the space in the upper-left
	is the closure space consisting of a single point.
	\begin{equation*}
		\begin{tikzcd}
			(*,c_*) \arrow[r] \arrow[d] & (E,*) \arrow[d,"p"]\\
		(X,*) \arrow[ur,"g",dashed] \arrow[r,"f"] & (B,*)
		\end{tikzcd}
	\end{equation*}
	Suppose, furthermore, that for
	every $b\in B$ there exists a neighborhood $U_b \subset B$ of $b$ and
neighborhoods $W_{\tilde{b}} \subset E$ of each $\tilde{b} \in p^{-1}(b)$
	such that
	\begin{enumerate}
		\item $W_{\tilde{b}} \cap W_{\tilde{b}'} = \emptyset$ for any
			$\tilde{b} \neq \tilde{b}'$, $\tilde{b},\tilde{b}' \in
			p^{-1}(b)$, and
		\item The restriction of $p$ to each $W_{\tilde{b}}$ is
			injective into $U_b$.
	\end{enumerate}
	Then the lift $g$ of $f$ is unique.
\end{proposition}

\begin{proof}
	Let $g,g':(X,c_X,*) \to (E,c_E,*)$ be lifts of $f$ with $g(*) = g'(*)$. Let $x \in X$, $\tilde{b} = g(x) \in E$ and
	$\tilde{b'}=g'(x) \in E$ so that $p(\tilde{b}) = p(\tilde{b}') = f(x)$, and consider the sets $W_{\tilde{b}}$ and
	$W_{\tilde{b}'}$. Since $W_{\tilde{b}}$ and
	$W_{\tilde{b}'}$ are neighborhoods of $\tilde{b}$ and 
	$\tilde{b}'$, respectively, it follows that $N = g^{-1}(W_{\tilde{b}}) \cap
	g'^{-1}(W_{\tilde{b}'})$ is a neighborhood of $x \in X$. If $g(x) =
	g'(x)$ then $g|N = g'|N$, since every point in
	$W_{\tilde{b}}$ and $W_{\tilde{b}'}$ is the pre-image of at most one point
	in $U_b$. Conversely, if $g(x) \neq g'(x)$, then $f(N) \cap g(N)
	\subset W_{\tilde{b}} \cap
	W_{\tilde{b}'} = \emptyset$.

	It follows that the sets $\{x \in X \mid g(x) = g'(x) \}$ and $\{x \in
	X \mid g(x) \neq g'(x) \}$ are both open and closed, and therefore they are
	either equal to $X$ or $\emptyset$, since $X$ is connected. Since $g(*)
	= g'(*)$, we have that $X = \{x \in X \mid g(x) = g'(x) \}$, and
	therefore $g = g'$ as
	desired.
\end{proof}

In the next propositions, we give severeal examples of covering maps.

\begin{proposition}
	\label{prop:S1 covering map}
	The map $p:(\R,c_r) \to (S^1,c_r)$, \[p(x) = x \mod 1 = [x]\] is a covering map iff $0\leq r < \frac{1}{3}$.
\end{proposition}
\begin{proof}
	Let $0 \leq r < \frac{1}{3}$. Let $x \in [0,1)$, so $p(x) = [x] \in S^1$.
	Define the set $U_x \coloneqq \{[y] \in S^1 \mid x-\frac{1}{3} < y <
	x+\frac{1}{3}\}$, and note that $U_x$ is a neighborhood of $[x]$. For
	every $n \in \Z$, define
	\[ U^n_x = \left(x + n - \tfrac{1}{3},x+n+\tfrac{1}{3}\right).\] 
	Then $(U^n_x,c_r)\cong (U_x,c_r)$ are homeomorphic as closure spaces,
	and, as a set, \[p^{-1}(U_x) = \coprod_{n \in \Z} U^n_x.\] Furthermore,
	the inclusion map $i:\left(\coprod_{n \in \Z} U^n_x,c_r\right)\to
	p^{-1}(U_x)$ is continuous, since $\left(\coprod_{n \in \Z}
	U^n_x,c_r\right)$ is a colimit in $\catname{Cl}$, and $i$ is the map
	guaranteed by the universal property of colimits. Since, for any $n\neq
	n'$, we have $d(U^n_x,U^{n'}_x) \geq \frac{1}{3}$, it follows that
	$i^{-1}$ is continuous as well, and therefore we have the closure space homeomorphism
	$(p^{-1}(U_x),c_r) \cong \left(\coprod_{n \in \Z} U^n_x,c_r\right)$.
	Since $\left(\coprod_{n \in \Z} U^n_x,c_r\right) \cong (\Z \times
	U_x,c_\Pi)$
	as closure spaces, where $\Z$ is given the discrete closure structure
	and $c_\Pi$ is the product closure structure, it follows that $p$ is a covering map.
	
	Now consider the case $r \geq \frac{1}{3}$. By Theorem \ref{thm:Connectedness}, a closure space $(X,c_X)$ is connected iff it is not the union of two nonempty disjoint open subsets, or, equivalently, that there is no proper subset of $X$ that is both open and closed. Now note that, with $r \geq \frac{1}{3}$, if $U$ is a neighborhood of a point $x \in S^1$, then $p^{-1}(U)$ is connected with respect to the subspace convergence structure. However, each sheet of $U \times p^{-1}(x)$ is both open and closed, and therefore $U \times p^{-1}(x)$ is not connected. Therefore, $U \times p^{-1}(x)$  and $p^{-1}(U)$ cannot be homeomorphic, and $p$ is not a covering map.
\end{proof}

\begin{definition}
	Let $p:(E,c_e) \to (B,c_B)$ and $f:(X,c_X) \to (B,c_B)$ be continuous
	maps. Define the set $X \times_B E \subset X\times E$ by 
\begin{equation*}
	X \times_B E \coloneqq \{(x,e) \in X \times E \mid f(x) = p(e)\}.
\end{equation*}
The closure structure on $X \times_B E$, denoted $c_{\times_B}$, is given by the subset closure
	structure.
\end{definition}

\begin{proposition}
	$(X \times_B E,c_{\times_B})$ is the pullback of the diagram
\begin{equation*}
\begin{tikzcd}
	& (E,c_E) \arrow[d,"p"]\\
	(X,c_X) \arrow[r,"f"] & (B,c_B)
\end{tikzcd}
\end{equation*}
\end{proposition}
\begin{proof}
	Let $(Y,c_Y)$ be a closure space, and suppose that $g:Y \to X$ and $h:Y
	\to E$ are continuous and make the diagram
\begin{equation*}
\begin{tikzcd}
	Y \arrow[r,"h"] \arrow[d,"g"]& E \arrow[d,"p"]\\
X \arrow[r,"f"] & B
\end{tikzcd}
\end{equation*}
commute. Then there is a map $\phi:Y \to X \times_B E$ of sets given by $\phi(y) =
(g(y),h(y))$. Since $g$ and $h$ are continuous, $\phi$ is continuous as a map
to the product space $(X\times E,c_{X\times E})$, i.e. as a map
$(Y,c_Y) \to (X\times E,c_{X\times E})$. However, since the range of $\phi$ is contained in $X
\times_B E$, it follows that, for any $Y \subset X$, we have
\begin{equation}
	\phi(c_Y(A)) \subset c_{X \times E}(\phi(A)) \cap (X \times_B E) =
	c_{X_B}(\phi(A)).
\end{equation}
Therefore, $\phi$ is continuous as a map from $(Y,c_Y) \to (X\times_B E,c_{\times_B})$, and the result
follows.
\end{proof}

\begin{theorem}
\label{thm:Pullback of covering map}
	The pullback of a covering map is a covering map.
\end{theorem}
\begin{proof}
Let $p:(E,c_E) \to (B,c_B)$ be a covering map, suppose the map $f:(X,c_X) \to
(B,c_B)$ is continuous, and let 
$(X \times_B E,c_{\times_B})$ be the pullback of the resulting diagram, and denote by
$q:(X\times_B E,c_{\times_B})\to (X,c_X)$ the resulting projection. Let $b = f(x)$ for some $x \in X$. Since $p$ is a covering map, 
there exists a neighborhood $U_b \subset B$ of $b$ such that 
$p^{-1}(U_b) \cong F \times U_b$. 
Define $V_x \coloneqq f^{-1}(U_b) \subset X$. Then 
\begin{align*}
	q^{-1}(V_x) &= \{(x,e) \in X \times E \mid f(x) = p(e) \in U_b \} \\
		    & = \{(x,e) \in X \times E \mid e \in p^{-1}(f(x)), f(x)
		    \in U_b \}\\
		    &= f^{-1}(U_b) \times F\\
		    &= V_x \times F.
\end{align*}
Therefore, $q:(X \times_B E,c_{\times_B}) \to (X,c_X)$ is a covering map.
\end{proof}

\begin{definition}
	For every $n \in \N$, let $(\Z_n,c_m)$ denote the closure space with
	$\Z_n = \Z/n\Z$, where the closure operator $c_m$ is given by
\begin{align*}
&c_m([k]) = \{[k-m], \dots, [k-1], [k], [k+1], \dots, [k+m]\} \\
&c_m(A) = \cup_{[k] \in A} c_m([k]).
\end{align*}
Furthermore, let $(\Z,c_m)$ denote the closure space consisting of $\Z$ with
the closure operator
\begin{align*}
&c_m(k) = \{k-m, \dots, k-1, k, k+1, \dots, k+m\} \\
&c_m(A) = \cup_{k \in A} c_m(k).
\end{align*} 
\end{definition}

\begin{proposition}
\label{prop:Pullback map} Let $\iota:(\Z_n,c_m) \to (S^1,c_{m/n})$ be the map
$\iota[k] = \left[\frac{k}{n}\right]$. Then $\iota$ is continuous, and
the map $q:(\Z,c_m) \to (\Z_n,c_m)$ given by $q(k) = [k]$ is the pullback of
	$p:(\R,c_{m/n}) \to (S^1,c_{m/n})$ along $\iota$.
\end{proposition}
\begin{proof}
	We first note that $\iota$ is well defined, since for any two
	representatives of $[k]$, we have $\iota[k] = \left[\frac{k}{n}\right]
	= \left[\frac{k}{n} + a\right] = \iota[k + an]$, for any $a\in \Z$. Next, since for any point $[k] \in \Z_n$, $\iota(c_m[k]) \subset
	c_{m/n}(\iota[k])$ by definition, it follows that $\iota$ is
	continuous.

	Let $(X,c_X)$ be a closure space, and suppose that $f:(X,c_X) \to
	(\R,c_{m/n})$ and $g:(X,c_X) \to (\Z_n,c_m)$ are continuous maps such
	that the diagram
\begin{equation*}
\begin{tikzcd}
	X \arrow[r,"f"] \arrow[d,"g"] & \R \arrow[d,"p"]\\
	\Z_n \arrow[r,"\iota"] & S^1
\end{tikzcd}
\end{equation*}
commutes.
	Define a map $\phi:(X,c_X) \to (\Z,c_m)$ by $n \cdot f(x)$. First, note that,
	since $p(f(x)) = \iota(g(x)) = [\frac{g(x)}{n}]$, where $g(x)$ is an
	integer, it follows that $f(x) = k + \frac{g(x)}{n}$ for some integer
	$k$. Therefore, $n\cdot f(x)$ is an integer for all $x\in X$, so $\phi$ is well-defined, and,
	furthermore, $f = \tilde{\iota}\circ \phi$ and $g = q \circ \phi$,
	where $\tilde{\iota}$ is the continuous map $\tilde{\iota}(k) =
	\frac{k}{n} \in \R$, and $q$ is as in the hypothesis.

	Furtherore, since $f$ is continuous, we have, for any $A\subset X$, that $f(c_X(A))
	\subseteq c_{m/n}(f(A))$. However, this implies that
	\[\phi(c_X(A))
	= n\cdot f(c_X(A)) \subseteq (n \cdot c_{m/n}(f(A))) \cap \Z  = c_m(n
\cdot f(A)) \cap \Z = c_m(\phi(A)),\] 
	and therefore $\phi$ is continuous. It
	follows that $(\Z,c_m)$ is the pullback of $p$ along $\iota$. By the
	unicity of pullbacks, we further have $(\Z,c_M) \cong (\Z \times_{S^1}
	\R,c_{\times_{S^1}})$.
\end{proof}

\begin{corollary}
	\label{prop:Zn covering map} The map $q:(\Z,c_m) \to (\Z_n,c_m)$ in
	Proposition \ref{prop:Pullback map} above is a covering map iff $0 \leq m < \frac{n}{3}$.
\end{corollary}

\begin{proof}
	By Proposition \ref{prop:Pullback map}, the map $q$ is the pullback of the
	covering map $p:(\R,c_{m/n}) \to (S^1,c_{m/n})$. By \ref{prop:S1
	covering map}, $p$ is a covering map iff $0\leq \frac{m}{n} <
	\frac{1}{3}$. By Theorem \ref{thm:Pullback of covering map}, $q$ is
	therefore a covering map if $0 \leq m < \frac{n}{3}$.

	Conversely, if $m \geq \frac{n}{3}$, the space $q^{-1}(U)$ is path
	connected for any neighborhood $U$ of any point $[k] \in \Z_n$, and
	therefore $q$ is not a covering map.
\end{proof}

The remainder of the calculation of the fundamental group of $(S^1,c_r)$ for
$0 < r < \frac{1}{3}$ now reduces to a generalization of the classical
calculation of the fundamental group of $S^1$. The calculation for the
fundamental group of $(\Z_n,c_m)$ will proceed similarly, although clearly this
case is specific to the closure space setting. We present the calculations here
as corollaries of the theorem that, given a covering map $p:E \to B$, where $E$
is a simply connected space, $\pi_1(X,*)$ is isomorphic to the automorphism
group of the covering. The proof is nearly verbatim the classical proof for
topological spaces, with a few delicate points regarding when to appeal to open
sets of a closure structure and when to use neighborhoods. We have adapted our
treatment from \cite{tom_Dieck_2008} and \cite{Spanier_1966}, and we include
several auxiliary results on covering spaces. We begin with the following
preparatory results.

\begin{lemma}
	\label{lem:Sheets are neighborhoods}
	Let $p:(E,c_E) \to (B,c_B)$ be a covering map, and let $U$ be a neighborhood of $b \in (B,c_B)$ such that $p^{-1}(U) \cong U \times p^{-1}(b)$. Let $x \in p^{-1}(b)$. Then the sheet $\tilde{U}_x \subset p^{-1}(U)$ which contains $x$ is a neighborhood of $x$ in $(E,c_E)$. 
\end{lemma}
\begin{proof}
	This follows immediately from Proposition \ref{prop:Neighborhood in
	neighborhood is neighborhood} with $X = E$, $U = p^{-1}(U)$, and $V =
	\tilde{U}_x$.
\end{proof}

\begin{proposition}
\label{prop:Uniqueness of liftings of coverings}
Let $p:(E,c_E) \to (B,c_B)$ be a covering map, and let \[g_0, g_1:(X,c_X) \to
(E,c_E)\] be liftings of a continuous map $f:(X,c_X) \to (B,c_B)$. Suppose $g_0(x)=g_1(x)$ for some $x \in X$. If $(X,c_X)$ is connected, then $g_0 = g_1$.
\end{proposition}
\begin{proof}
	Since $p$ is a covering map, the hypotheses of Proposition
	\ref{prop:Unique liftings intermediate} are satisfied. The result
	follows.
\end{proof}

\begin{lemma}
\label{lem:Product open sets}
Let $\mathcal{U}$ be an an interior cover of $X \times I$ with the product closure structure. For each $x\in X$ there exists a neighborhood $V(x)$ of $x \in X$ and $n = n(x) \in \mathbb{N}$ such that, for $0 \leq i < n$, the set $V(x) \times \left[\frac{i}{n},\frac{(i+1)}{n}\right]$ is contained in some member of $\mathcal{U}$.
\end{lemma}
\begin{proof}
	By Theorem \ref{thm:Product local base and local subbase}, if $\mathcal{U}$ is a local subbase for $x\in X$, and $\mathcal{V}$ is a local subbase for $t \in [0,1]$, then the family $\{\pi_1^{-1}(U) \mid U\in \mathcal{U}\} \cup \{\pi_2^{-1}(V) \mid V \in \mathcal{V}\}$ is a local subbase for $(x,t) \in X\times I$, where the $\pi_i$, $i \in \{1,2\}$ are the projections onto $X$ and $I$, respectively. It follows that every neighborhood $N_t \subset X \times I$ of a point $(x,t) \in X \times I$ contains a neighborhood of $(x,t)$ of the form $V_{(x,t)} = U_t(x) \times (t_1,t_2) \subset X \times I$, where $U_t(x)$ is a neighborhood of $x \in X$ and $(t_1,t_2)$ is a neighborhood of $t \in [0,1]$. In particular, for every $(x,t)$, there is a neighborhood of the form of $V_{(x,t)}$ contained in any set of $\mathcal{N}$ which is a neighborhood of $(x,t)$. Since $I$ is compact, a finite number of such neighborhoods $V_{(x,t)}$ cover $\{x\} \times [0,1]$. Let $\{t_0, \dots, t_k\}$ be the points in $I$ such that the family $\{V_{(x,t_i)}\}_{i=0}^k$ covers $\{x\} \times I$, and let $\lambda$ be the Lebesgue number of the cover  $\{\pi_2(V_{(x,t_i)})\}_{i=0}^k$ of $I$ given by the projection onto $I$. Choosing $V(x) \coloneqq \cap_{i=1}^k U_{t_i}(x)$ and $n > \frac{1}{\lambda}$, the lemma follows.   
\end{proof}

We now give the definition of a fibration in $\catname{Cl}$ for use in what follows.
\begin{definition}
	We say that a map between closure spaces $p:(E,c_E) \to (B,c_B)$ has the \emph{homotopy lifting property} for the space $(X,c_X)$, if, for any homotopy $H:(X,c_X) \times I\to (X,c_X)$, and any map $f:(X,c_X) \to (E,c_E)$ such that $pf(x) = h(x,0)$, there exists a map $\tilde{H}:(X,c_X) \times I \to (E,c_E)$ such that $p\tilde{H} = H$. That is, for any commutative solid arrow diagram
\begin{equation*}
\begin{tikzcd}
X \ar [r,"f"] \ar [d,"i_0"] & E \ar [d,"p"] \\
X \times I \ar [r,"h"] \ar [dashed,ur,"H"] & B
\end{tikzcd}
\end{equation*}
there exists a dashed arrow making the complete diagram commute. The map $p$ is a \emph{fibration} iff it has the homotopy lifting property for all closure spaces $(X,c_X)$.
\end{definition}

\begin{lemma}
\label{lem:Projection is a fibration}
A projection $p:(B \times F,c_{B\times F}) \to (B,c_B)$ is a fibration.
\end{lemma}
\begin{proof}
Let $h:(X \times I, c_{X \times I}) \to (B,c_B)$ be a continuous map and let $h_0:(X,c_X) \to (B \times F,c_{B\times F})$ satisfy $ph_0(x) = h(x,0)$. Define $H:(X \times I, c_{X \times I}) \to (B \times F,c_{B\times F})$ by $H(x,t) = (h(x,t),\pi_F(h_0(x)))$, where $\pi_F:(B \times F,c_{B\times F}) \to (F,c_F)$ is the projection. Then $pH = h$ and $H$ is a lifting of $h$. Since $h$ and $X$ are arbitrary, $p$ is a fibration. 
\end{proof}
\begin{proposition}
\label{prop:Covering is a fibration}
A covering map $p:(E,c_E) \to (B,c_B)$ is a fibration.
\end{proposition}
\begin{proof}
Let $h:(X,c_X) \times I \to (B,c_B)$ and $a:(X,c_X) \to (E,c_E)$ be a homotopy and initial condition, respectively, i.e. $pa(x) = hi_0(x)$, where $i_0:X \to X \times I$ is the map $x \mapsto (x,0)$.

Since $p$ is a covering map, there is a neighborhood system $\mathcal{U}_B$ on $B$ such that for each $U \in \mathcal{U}_B$, $p^{-1}(U) = U \times F_b$ for some $b\in U$. Since $h$ is continuous, $\mathcal{V} \coloneqq h^{-1}(\mathcal{U})$ is a neighborhood system on $X \times I$.  Consider a point $x \in X$, and let $V(x) \subset X$ and $n_x \in \mathcal{N}$ be the neighborhood of $X$ and the natural number guaranteed by Lemma \ref{lem:Product open sets}, respectively, so that $V(x) \times \left[\frac{i}{n_x},\frac{(i+1)}{n_x}
\right] \subset V$ for some $V \in \mathcal{V}$.
 
We proceed by induction on $i$. First consider the case $i=0$. Since $p:p^{-1}(U) \to U$ is a projection, it is a fibration by \ref{lem:Projection is a fibration}, and therefore $h:V(x) \times \left[0,\frac{1}{n_x}\right] \to U$ has a lifting $H:V(x) \times \left[0,\frac{1}{n_x}\right] \to E$ with initial condition $a|_{V(x)}$. 

Now suppose there exists a lifting of $h|_{V(x) \times \left[0,\frac{i}{n_x}\right]} = h_i:V(x) \times 
\left[0,\frac{i}{n_x}\right] \to (B,c_B)$. Given such a lifting $H_i:V(x) \times \left[0,\frac{i}{n_x}\right] = (E,c_E)$, this 
same argument gives a lifting $H'_{i+1}:V(x) \times \left[\frac{i}{n_x},\frac{i+1}{n_x}\right]\to (E,c_E)$ of $h|{V(x) \times 
\left[\frac{i}{n_x},\frac{i+1}{n_x}\right]}$ with initial condition $H|_{V(x) \times \left[\frac{i}{n_x}\right]}$. Since the lifting of $h(y \times I)$ is unique, The maps $H_{i}$ and $H'_{i+1}$ agree on $V(x) \times \frac{i}{n_x}$, the map $H_{i+1}:V(x) \times \left[0,\frac{i+1}{n_x}\right] \to (E,c_E)$ given by
\begin{equation*}
H_{i+1}(y,t) = \begin{cases}
H_{i}(y,t) & t \in \left[0,\frac{i}{n_x}\right]\\
H'_{i+1}(y,t) & t \in \left[\frac{i}{n_x},\frac{i+1}{n_x}\right]
\end{cases}
\end{equation*}
is continuous by \ref{cor:Cont function from restriction to closed subspaces} as a map from the subspace $V(x) \times [0,\frac{i+1}{n_x}] \subset (X,c_X)$ with the induced closure operator. $I$ is connected, Proposition \ref{prop:Uniqueness of liftings of coverings} shows that the lifting of $h|x \times I$ to $E$ is unique for any $x$, and therefore 
the liftings on each $V(x) \times I$ combine to a well-defined lifting $H:X \times I \to E$ of $h$. Since the $V(x)\times I$ form an interior cover of $(X,c_X) \times I$, we have that $H$ is continuous by \ref{cor:Cont function from restriction to interior cover}. Since $(X,c_X)$ was arbitrary, $p$ is a fibration.
\end{proof}

The remainder of the results required for our calculations are primarily algebraic, and follow verbatim the discussion in \cite{tom_Dieck_2008}, Section 3.2, replacing topological spaces with closure spaces. We give a full exposition here for convenience, as well as to establish notation and definitions for the computations which follow.

\begin{lemma}
	\label{lem:HLP alternate}
	Let $p:(E,c_E) \to (B,c_B)$ have the homotopy lifting property for the cube $I^n$. Then for each commutative diagram
	\begin{equation*}
	\begin{tikzcd}
	I^n \times 0 \cup \partial I^n \times I \ar [d,"i"] \ar [r,"a"] & E \ar [d,"p"]\\
	I^n \times I \ar [ur, dashed] \ar [r,"h"] & B
	\end{tikzcd}
	\end{equation*}
	where $i$ is the inclusion, there exists $H:I^n \times I \to (E,c_E)$ with $Hi = a$ and $pH = h$.	
\end{lemma}
\begin{proof} We first recall that there exists a homeomorphism \[\phi:(I^n \times I, I^n \times \{0\}) \to (I^n \times I, (I^n \times 0) \cup (\partial I^n \times I)).\] Using this homeomorphism, we may transform the above lifting problem into a homotopy lifting problem for the cube $I^n$. Since $p$ has the homotopy lifting property for cubes, a function $H: I^n \times I$ exists with the desired properties.
\end{proof}

\begin{definition}
Let $(A,a),(B,b)$, and $(C,c)$ be pointed sets. We say that a sequence $(A,a) \xrightarrow{f} (B,b) \xrightarrow{g} (C,c)$ is \emph{exact} at $(B,b)$ iff $im (f) = \ker g$, where $\ker g \coloneqq g^{-1}(c)$. For a group, we take the basepoint to be the identity, and for the set $\pi_0(X)$ for some closure space $(X,c)$, we write $\pi_0(X,a)$, to indicate that the basepoint is $[a]$.
\end{definition}

\begin{theorem}
\label{thm:Fibration exact sequence}
Let $p:(E,c_E) \to (B,c_B)$ be a map with the homotopy lifting property for a point and for $I$, and write $F_b = p^{-1}(b)$. Then there exists a map $\partial_x:\pi_1(B,p(x)) \to \pi_0(F_b,x)$ making the sequence
\begin{equation*}
\begin{tikzcd}[column sep=15pt]
\pi_1(F_b,x) \ar [r,"i_*"] & \pi_1(E,x) \ar [r,"p_*"] & \pi_1(B,p(x)) \ar [r,"\partial_x"] & \pi_0(F_b,x) \ar [r,"i_*"] & \pi_0(E,x) \ar [r,"p_*"] & \pi_0(B,b)
\end{tikzcd} 
\end{equation*}
exact. 

Furthermore, the preimages of elements under $\partial_x$, i.e. $\partial_x^{-1}[x']$ for $[x'] \in \pi(F_b,x)$, are the left cosets of $\pi_1(B,b)$ with respect to $p_*\pi_1(E,x)$. Finally, the preimages of $\pi_0(i)=i_*:\pi_0(F_b,x)\to \pi_0(E,x)$ are the orbits of the $\pi_1(B,b)$-action on $\pi_0(F_b,x)$.
\end{theorem}

\begin{proof}
Let $v:I \to B$ be a path from $v(0) = b \in B$ to $v(1) = c \in B$. For each such path 
we will define a map $v_{\#}:\pi_0(F_b) \to \pi_0(F_c)$. Let $x \in F_b$. Since $p$ has the homotopy lifting property for $I$, we may find a lifting $V:I \to 
E$ of $v$ with $V(0) = x$. Define $v_{\#}[x] = [V(1)]$. 

We claim that $v_{\#}$ is  well-defined and depends only on the class $[v_*] \in \Pi(B)$. Let $x'\in F_b$ with 
$[x']=[x] \in F_b$. Then there is a path $\gamma:I \to F_b$ which connects $x$ to $x'$. 
Let $v':I \to B$ be a path from $b$ to $c$ with $v \simeq v'\text{ rel }\partial I$, and let $h:I \times I \to B$ be the homotopy from $v$ to $v'$. Finally, let $V, V': I \to E$ be liftings of $v$ and $v'$ to $E$ with initial points $\gamma(0)$ and $\gamma(1)$, 
respectively. The maps $\gamma, V$, and $V'$ give a continuous map $a:I \times \partial I \cup 0 \times I \to E$, defined by $a(s,0) = V(s)$, $a(s,1) = V'(s)$, $a(0,t) = \gamma(t)$. We also have that $p\circ a = h \circ i$, where $i:I \times \partial I \cup 0 \times I  \to I \times I$ is the inclusion. Then Lemma \ref{lem:HLP alternate} gives a lifting $H$ of $h$ with initial condition $a$. Since $p\circ H(1,t) = c$, the homotopy $H$ defines a path $H(1,t)$ from $V(1) \to V'(1)$ in $F_c$, and therefore $[V(1)] = [V'(1)] \in \pi_0(F_c)$. This shows that the map $v_{\#}$ is well defined and depends only on the class $[v] \in \Pi(X)$, proving the claim. Additionally, it follows from the definitions that $w_{\#}v_{\#} = (v \star w)_{\#}$

Note that map $v_{\#}$ gives a right action of $\pi_1(B,b)$ on $\pi_0(F_b)$ by $[x]\cdot [v] \mapsto v_{\#}[x]$. We use this to define a map $\partial_x: \pi_1(B,b) \to \pi_0(F_b,x)$ by
\begin{align*}
\partial_x[v] = [x] \cdot [v] = v_{\#}[x].
\end{align*}

With this definition, we turn our attention to the exactness of the sequence in the theorem. By the definitions of the maps, we see that the composition of any two of them is equal to the basepoint. We now check that the kernel of every map is contained in the image of the previous one.

Let $[u] \in \ker p_* \subset \pi_1(E,x)$ and let $h: I \times I \to B$ be a null-homotopy of $p\circ u$. Let $a:I \times 0 \cup \partial I \times I \to E$ be the map $a(s,0) = u(s)$ and $a(\epsilon,t) = x$, $\epsilon = 0,1$, and consider the lifing problem for $h$ with initial condition $a$. The lifting $H$ of $h$ is a homotopy of loops beginning at $u$ and ending at a loop $u'$ with $p\circ u' = b$, which implies that $[u']$ is in the image of $i_*$.

Let $\delta_x[v] = [x]$. Then there exists a lifting $V$ of $v$ with $V(1) \in [x]$. Choose a path $\gamma:I \to F_b$ from $V(1) \to x$. Since $V(0)= x$ by the definition of $\delta_x$, it follows that $V \star \gamma$ is a loop in $E$, and $V$ satisfies $p_*[V \times \gamma] = [v]$, since $p_*[\gamma] = [b]$ is constant. Therefore $[v]$ is in the image of $p_*$. 

Let $i_*[y] = [x] \in \pi_0(E,x)$. Then there exists a path $\gamma:I \to E$ with $\gamma(0) = x$ and $\gamma(1) = y$. The projection $p\circ \gamma$ is a loop in $B$ based at $b$, and $\partial_x[p\circ \gamma] = [\gamma(1)] = [y]$. Therefore $[y]$ is in the image of $\partial_x$.

Let $p_*[y] = [b] \in \pi_0(B,b)$. Let $\gamma:I \to B$ be a path from $p(y)$ to $b$. Let $\Gamma:I \to E$ be a lifting of $\gamma$ with $\Gamma(0) = y$. Then $\Gamma(1) \in p^{-1}(b) = F_b$, and therefore $[y] = i_*[\Gamma(1)] \in im(i_*)$.

We have now shown that the sequence in the theorem is exact, and we proceed to prove the final statements.
Let $\partial_x[u] = \partial_x[v]$, where $[u],[v] \in \pi_1(B,b)$, and choose liftings $U$,$V$ of $u$ and $v$, respectively, with $U(0) = x = V(0)$. Since $[U(1)] = \partial_x[u] = \partial_x[v] = [V(1)] \in \pi_0(F_b,x)$, there exists a path $\gamma:I \to F_b$ from $U(1)$ to $V(1)$. Then $U \star \gamma \star V^{-1}$ is a loop in $E$, and
\begin{align*}
[p_*(U \star \gamma \star V^{-1}) \star v] &= [p_*(U) \star p_*(\gamma) \star p_*(V^{-1})][v]\\
& = [u][b][v^{-1}][v]\\
& = [u],
\end{align*}  
where $b$ above is the constant map $b:I \to b \in B$. Therefore, $[u], [v]\in \pi_1(B,b)$ are in the same left coset. 

Conversely, if $[u]$ and $[v]$ are in the same coset, then there exists a $[\gamma] \in \pi_1(E,x)$ such that $p_*[\gamma][v] = [u] \in \pi_0(F_b,b)$. This implies that $(p \circ \gamma) \star v \simeq u$, so let $h: I \times I \to B$ the a homotopy between $(p \circ \gamma) \star v$ and $u$. Let $U,V$ be liftings of $u$ and $(p \circ \gamma) \star v$, respectively, with $U(0) = x = V(0)$. Then there exists a map $H: I \times I \to E$ which lifts $h$ and has initial data $a:I \times 0 \cup \partial I \times I$ given by $a(s,0) = \gamma(0) = x$, $a(0,t) = \gamma\star V(t)$, and $a(1,t) = U(t)$. Then $H(1,t)$ is a path from $U(1)$ to $V(1)$ contained in $F_b$, and therefore $\delta_x[u] = \delta_x[v]$, as desired. 

Now suppose that $i_*[a] = i_*[b] \in \pi_0(E,x)$. Then there exists a path $\gamma:I \to E$ from $a$ to $b$, and $p\circ \gamma$ is a loop in $B$ with basepoint $b$. Then, by definition of the right action of $\pi_1(B,b)$ on $\pi_0(F_b,b)$, we have $[a]\cdot [p\circ \gamma] = [b]$, and therefore $[a]$ and $[b]$ are in the same orbit. 

Conversely, suppose that $[a],[b] \in \pi_0(F_b,x)$ are in the same orbit of $\pi_1(B,b)$. Then $[a]\cdot[\gamma] = [b]$ for some $[\gamma] \in \pi_1(B,b)$, and therefore there exists a lifting $\Gamma$ of $\gamma$ to $E$ such that $[a]\cdot [\gamma] = [\Gamma(1)] = [b]$ and $\Gamma(0) = a$. Since there exists a path in $F_b$ between $\Gamma(1)$ and $b$, it follows that $i_*[a] = i_*[b]$, and the proof is complete.
\end{proof}

Specializing to the case where $(B,c_B)$ is path-connected and $p$ is a covering, we have

\begin{corollary}
	\label{cor:Exact sequence covering}
	Let $p:(E,c_E) \to (B,c_B)$ be a covering, with $B$ path-connected. Then the sequence
	\begin{equation*}
	\begin{tikzcd}[column sep=20pt]
	1 \ar[r] & \pi_1(E,x) \ar [r,"p_*"] & \pi_1(B,p(x)) \ar [r,"\partial_x"]& \pi_0(F_b,x) \ar [r,"i_*"] & \pi_0(E,x) \ar [r] & 1
		\end{tikzcd}
	\end{equation*}
	is exact and $i_*$ is surjective. Therefore, E is path connected iff $\pi_1(B,b)$ acts transitively on $F_b$. The isotropy group of $x \in F_b$ is the image of $p_*:\pi_1(E,x) \to \pi_1(B,b)$.
\end{corollary}
\begin{proof}
	Since $F_b$ is discrete and $B$ is path-connected, $\pi_1(F_b,b) = 1$ and $\pi_0(B,b) = 1$. The result now follows from Theorem \ref{thm:Fibration exact sequence}. 
\end{proof}

We now show that, when there is a group action on the total space $E$ which leaves $p:E \to B$ invariant, we may replace the above sequence with an exact sequence of groups. 

\begin{lemma}
	\label{lem:Fibers of coverings have constant paths}
	Let $p:(E,c_E) \to (B,c_B)$ be a covering map. Then any path $\gamma:I \to p^{-1}(x)$ is constant for any $x\in B$.
\end{lemma}
\begin{proof}
	Let $\gamma:I \to p^{-1}(x)$ be a path, and let $\gamma':I \to p^{-1}(x)$ be the constant map $\gamma'(t) = \gamma(0)$. Since they both cover the map $f:I \to B$, $f(t) = x$ and have the same initial point $\gamma(0) = \gamma'(0)$, we must have $\gamma'(t) = \gamma(t) = \gamma(0)$ by Proposition \ref{prop:Uniqueness of liftings of coverings}.
\end{proof}

\begin{proposition}
	\label{prop:Liftings of homotopies in coverings}
	Let $p:(E,c_E) \to (B,c_B)$ be a covering map. Let $\gamma:I \to (B,c_B)$ be a path with starting point $\gamma(0) = p(e)$ for some $e\in E$. Then there exists a unique lifting $\Gamma:I \to (E,c_E)$ of $w$ with $\Gamma(0) = e$. Furthermore, if two paths $\gamma_0, \gamma_1:I \to (E,c_E)$ start at the same point $\gamma_0(0) = \gamma_1(0)$ and $p\circ \gamma_0 \simeq p \circ \gamma_1 \text{ rel } \partial I$, then $\gamma_0 \simeq \gamma_1 \text{ rel } \partial I$. 
\end{proposition}

\begin{proof}
	The lifting exists by Proposition \ref{prop:Covering is a fibration}, and it is unique by Proposition \ref{prop:Uniqueness of liftings of coverings}. Now let $h:I \times I \to B$ be a homotopy between paths $p\circ \gamma_0, p\circ \gamma_1:I \to E$. By the homotopy lifting property, there exists a map $H:I \times I \to E$ such that $H(s,0) = \gamma_0(s)$. Observe that $H(0,t) \in p^{-1}(\gamma_0(0)$ and $H(1,t) \in p^{-1}(\gamma_1(1))$. By Lemma \ref{lem:Fibers of coverings have constant paths} $H(0,t)$ and $H(1,t)$ are constant. Therefore, $H$ is a homotopy rel $\partial I$ from $\gamma_0$ to a path $\gamma':I \to E$ with $\gamma'(0) = \gamma_1(0)$ and which lifts $\gamma_1$. By Proposition \ref{prop:Uniqueness of liftings of coverings}, $\gamma' = \gamma_1$.
\end{proof}

\begin{definition}
	Let $p:(E,c_e) \to (B,c_B)$ be a covering. A \emph{right $G$-action} $E \times G \to E$, $(e,g) \mapsto eg$ is said to be \emph{totally discontinuous} iff every $x \in E$ has a neighborhood $U$ with $Ug \cap U = \emptyset$ and $Ug \cup U$ is not connected for all $g \in G$. We call $p$ a \emph{G-principal covering} if $G$ is totally discontinuous and $p(xg) = p(x)$ for all $g \in G$ and such that the induced action on each fiber $F_b = p^{-1}(b)$ is transitive.
\end{definition}

\begin{theorem}
	\label{thm:G-principal covering}
	Let $p:(E,c_E) \to (B,c_B)$ be a right $G$-principal covering with path-connected total space $(E,c_E)$. Then the sequence of groups and homomorphisms
	\begin{equation*}
	\begin{tikzcd}
	1 \ar[r] & \pi_1(E,x) \ar [r,"p_*"] & \pi_1(B,p(x)) \ar [r,"\Gamma_x"]& G \ar [r] & 1 
	\end{tikzcd}
	\end{equation*}
	is exact (for each $x \in E$). The image of $p_*$ is a normal subgroup. The space $E$ is simply connected iff $\delta_x$ is an isomorphism. Thus, if $E$ is simply connected, $G$ is isomorphic to the fundamental group of $B$.
\end{theorem}

\begin{proof}
	We first recall the right action of $\pi_1(B,b)$ on $(F_b,x)$ given by $x\cdot [v] = v_{\#}(x) = V(1)$, where $V$ is a map which lifts $v$ and has initial point $x$. By Proposition \ref{prop:Liftings of homotopies in coverings}, the point $V(1)$ only depends on the homotopy class of $[v]$. Now let $g \in G$, and observe that $xg\cdot[v] = W(1)$, where $W:I \to E$ is a lift of $v$ which begins at $xg \in E$. We also have $(x\cdot [v])g = V'(1)g$, where $V':I\to E$ is a lift of $v$ which begins at $x$. It follows that $V'g$ is a lift of $v$ which begins at $gx$, and by Proposition \ref{prop:Uniqueness of liftings of coverings}, the uniqueness of liftings for covering maps, we have that $W = V'g$. Therefore, the actions of $G$ and $\pi_1(B,b)$ on $F_b$ commute. 
	
	Fix $x \in E$. Since $G$ is free and transitive on $F_b$, for every $[v]  \in \pi_1(B,b)$ and $x \in F_b$ there is a unique $\Gamma_x([v])\in G$ such that $x\cdot[v] = x\cdot \Gamma_x^{-1}([v])$. This defines a map $\Gamma_x:\pi_1(B,b) \to G$, $[v] \mapsto g$ such that $x[v] = x \Gamma_x([v])^{-1} = xg^{-1}$. To see that it is a homomorphism, we compute
	\begin{align*}
	x\cdot \Gamma([v][w])^{-1} = x \cdot [v][w]  = x\Gamma([v])^{-1}[w] = x [w] \Gamma([v])^{-1} = x \Gamma([w])^{-1}\Gamma([v])^{-1}
	\end{align*}
	Define a map $\rho_x:G \to F_b$ by $\rho_x(g) = x\cdot g^{-1}$. Since $G$ acts freely and transitively, $\rho$ is a bijection. Furthermore $\rho_x \Gamma_x = \partial_x$, i.e $\rho_x \Gamma_x [v] = x \cdot \Gamma_x([v])^{-1} = x \cdot [v] = \partial_x [v]$ for all $[v] \in \pi_1(B,b)$. 
	
	It now follows from Corollary \ref{cor:Exact sequence covering} that the sequence in the conclusion of the Theorem is exact, taking $\Gamma_x = \rho_x^{-1}\circ \partial_x$
\end{proof}

For the computation of the fundamental froup of $(S^1,c_r)$, we first need to establish that $\pi_1(\R,c_r) = 1$. As in the topological case, $\R$ endowed with these closure structures is contractible. 

\begin{lemma}
	\label{lem:R cr contractible}
	The space $(\R,c_r)$ is contractible, for any $r \geq 0$.
\end{lemma}
\begin{proof}
Observe first that, for any $x \in \R$, the collection of sets $\{(x - r -
\epsilon,x+r+\epsilon)\}$ forms a local base at $x$, Define $H:(\R,c_r) \times
I \to (\R,c_r)$ by $H(s,t) = st$. Then $H(s,0) = 0$ and $H(s,1) = s =
Id_{\R}(s)$. We claim that $H$ is continuous. To see this, fix $s \in \R$, $t
\in [0,1]$, and consider the interval $A \coloneqq (st-r-\epsilon,st + r + \epsilon)$. Then for $\alpha, \delta>0$, the set
\begin{equation*}
(s-r-\delta,s+r+\delta) \times (t-\alpha,t+\alpha)
\end{equation*}
is in $H^{-1}(A)$ iff 
\begin{align*}
((s-r-\delta)(t-\alpha),(s+r+\delta)(t+\alpha)) \subset (st - r - \epsilon, st + r + \epsilon)
\end{align*}
We compute
\begin{align*}
(s + r + \delta)(t + \beta) &= st + rt + \delta t + (s + r + \delta)\beta\\
&< st + r + \delta t +(s + r + \delta)\beta. 
\end{align*} 
Since $st + r < st + r + \epsilon$, there exist $\delta',\beta >0$ such that 
\[st + r + \delta' t +(s + r + \delta)\beta < st + r + \epsilon \]
Similarly, we may choose $\alpha$ and $\delta$ such that
\begin{align*}
st - r - \epsilon & < st - r - \delta t - (s-r-\delta)\alpha\\
& < st - rt - \delta t - (s - r - \delta)(\alpha)\\
& < (s-r-\delta)(t - \alpha)
\end{align*}
Choosing $\alpha, \delta>0$ to satisfy both of these constraints, we conclude
that $H$ is continuous at $(s,t)$. (Note that this argument also covers $t\in
\{0,1\}$.) Since $(s,t) \in \R \times I$ was abitrary, $H$ is therefore continuous. Since $H$ is a homotopy from $(\R,c_r)$ to the constant map, we conclude that $(\R,c_r)$ is contractible.
\end{proof}

\begin{theorem}
	\label{thm:Fundamental group of the circle 1}
	For $0 \leq r < \frac{1}{3}$, $\pi_1(S^1,c_r) \cong \Z$.
\end{theorem}
\begin{proof}
	For $r< \frac{1}{3}$, the map $p:(\R,c_r) \to (S^1,c_r)$ is a covering map by Proposition \ref{prop:S1 covering map}, and the map $\R \times \Z \to \R, (n,t) \mapsto t+n$ defines a totally discontinuous, right $\Z$-action $\R$ such that $p(t\cdot n) = p(t)$ for any $n \in \Z$. Since the induced action on each fiber is transitive, $p$ is a right $\Z$-principal covering map. Since $(\R,c_r)$ is contractible, the conclusion follows from Theorem \ref{thm:G-principal covering}. 
\end{proof}

We now compute $\pi_1(\Z_n,c_m)$ for $1 \leq m < \frac{n}{3}$. The proof is similar to the one above for $(S^1, c_r)$. We begin, as above, by proving the contractibility of the total spaces of the maps $q:(\Z,c_m) \to (\Z_n,c_m)$.

\begin{lemma}
	\label{lem:Z cm contractible}
	For any $m\geq 1$, $(\Z,c_m)$ is contractible.
\end{lemma}
\begin{proof}
	We first remark that a retract of a contractible space is contractible.
	Indeed, if the inclusion $\iota:(U,c_U) \to (X,c_X)$ has a left inverse
	$r:(X,c_X) \to (U,c_U)$ and $(X,c_X)$ is contractible, then the
	composition
	\begin{equation*}
	\begin{tikzcd}
		U \times I \ar[r,"\iota \times id"] & X \times I \ar[r,"H"] & X
		\ar[r,"r"] & U
	\end{tikzcd}
	\end{equation*}
	is a contraction of $U$ if $H$ is a contraction of $X$.

	Since $(\R,c_m)$ is contractible by Lemma \ref{lem:R cr contractible}, it
	remains to construct a retract $r:(\R,c_m) \to (\Z,c_m)$. Consider a
	strictly increasing sequence $\{a_k\}_{k\in \Z}$ such that each $a_k \in
	(0,1)$, i.e.
	\begin{equation*}
		0 < \cdots < a_{-2} < a_{-1} < a_0 < a_1 < a_2 < \cdots < 1
	\end{equation*}
	We define $r:(\R,c_m) \to (\Z,c_m)$ by
	\begin{equation*}
		r(x) = k, \text{ for } x\in [k-1+a_{k-1},k+a_k).
	\end{equation*}
	We will show that $r$ is continuous.

	Suppose $r(x) = k$. Then $x \in [k-1+a_{k-1},k+a_k)$, and choosing
	\[
	\delta < \min \{a_{k+m}-a_k,a_{k-1} - a_{k-1-m}\}
	\]
	we have that the set
	\begin{equation*}
		V \coloneqq\; [k-1+a_{k-1} - m - \delta, k+a_k + m + \delta)\\
	\end{equation*}
	is therefore a neighborhood of $x \in (\R,c_m)$. Observe that, by
	construction,
	\begin{equation*}
		V \subset \;[k-m-1+a_{k-m-1},k+m+a_{k+m}),
	\end{equation*}
	and it follows that 
	\begin{equation*}
		r(V) \subset \{k-m,k-m+1,\dots,k+m-1,k+m\}.
	\end{equation*}
	Since every neighborhood of $k$ in
	$(\Z,c_m)$ must contain $\{k-m,\dots,k+m\}$, Theorem
	\ref{thm:Neighborhood continuity} implies that $r$ is continuous, and therefore $(\Z,c_m)$ is
	contractible.
\end{proof}

\begin{corollary}
		\label{cor:Fundamental group of Zn}
		For $1 \leq m < \frac{n}{3}$, $\pi_1(\Z_n,c_m) \cong \Z$.
\end{corollary}
\begin{proof}
	For $1 \leq m < \frac{n}{3}$, the map $p:(\Z,c_m) \to (Z_n,c_m)$ is a covering map by Proposition \ref{prop:Zn covering map}, and the map $\Z \times \Z \to \Z, (k,j) \mapsto k+j$ defines a totally discontinuous, right $\Z$-action such that $p(k\cdot j) = p(k+j)=p(k)$ for any $n \in \Z$. Since the induced action on each fiber is also transitive, $p$ is a right $\Z$-principal covering map. Since $(\Z,c_m)$ is contractible, the conclusion follows from Theorem \ref{thm:G-principal covering}. 
\end{proof}

The above results show that $\pi_1(S^1,\tau) \cong \pi_1(Z_n,c_m)$ for $1 \leq
m < \frac{n}{3}$, but they do not give a map between the spaces which induces the isomorphism. The next result produces such a map. In particular, we show here the `nearest neighbor' maps, i.e. maps  $(S^1, \tau)$ and $(\Z_n,c_m)$ of the form found in Equation \ref{eq:Nearest neighbor map} of Section \ref{sec:Introduction}, induce an isomorphism on the fundamental groups.

\begin{theorem}
	Let $m,n\in \N$ satisfy $1\leq m < \frac{n}{3}$, and define $f:(S^1,\tau) \to (\Z_n,c_m)$ by
	\begin{equation*}
	f(t) = i, \text{   } t \in \left[\frac{2i-1}{2n},\frac{2i+1}{2n}\right),
	\end{equation*}
	where the interval is understood mod $1$. Then $f_*: \pi_1(S^1,\tau) \to \pi_1(\Z_n,c_m)$ is an isomorphism.
\end{theorem}
\begin{proof}
	Since $1\leq m < \frac{n}{3}$, the map $q:(\Z,c_m) \to (\Z_n,c_m)$, $q(k) = k \mod m$ is a covering map by Proposition \ref{prop:Zn covering map}, as is $p:(\R,\tau) \to (S^1,\tau)$, $p(x) = x \mod 1$ by Proposition \ref{prop:S1 covering map}. 
	We now consider the diagram
	\begin{equation*}
	\begin{tikzcd}
	(\Z,0) \arrow{d}{\tilde{f}|_{p^{-1}(0)} = \times n} \arrow{r}{i} &(\R,0) \arrow{r}{p}\arrow{d}{\tilde{f}} &(S^1,0) \arrow{d}{f}\\
	(n\Z,0)\arrow{r}{i} &(\Z,0) \arrow{r}{q}&(\Z_n,0) 
 	\end{tikzcd}
	\end{equation*}
	This map is commutative. Since $(\Z,c_m)$ and $(\R,c_{m/n})$ are
	contractible by Lemmas \ref{lem:R cr contractible} and \ref{lem:Z cm contractible}, it follows by Corollary \ref{cor:Exact sequence covering} that the maps in the above diagram induce the diagram
	\begin{equation*}
	\begin{tikzcd}
	1 \arrow{r} \arrow{d} & 1 \arrow{r}\arrow{d} & \pi_1(S^1,0) \arrow{r}{\partial_0} \arrow{d}{f_*}&\pi_0(\Z,0)\arrow{r} \arrow{d}{(\times n)_*}&1 \arrow{d}\\
	1 \arrow{r} & 1 \arrow{r} &\pi_1(\Z_n,0) \arrow{r}{\partial_0} &\pi_0(n\Z) \arrow{r} &1
	\end{tikzcd}
	\end{equation*}
	We claim that the non-trivial square of the diagram above is commutative. Let $[\gamma] \in \pi_1(S^1,0)$ Then $\partial_0 [\gamma] = [\Gamma(1)]$, where $\Gamma:I \to \R$ is a lift of $\gamma$ with $\Gamma(0) = 0$. Then $\tilde{f} \partial_0[\gamma] = [\tilde{f}\circ \Gamma(1)]$. Conversely, $\partial_0 f_*[\gamma] = [V(1)]$, where $V:I \to \Z$ is a lift of $f\circ \gamma$ with $V(0) = 0$. However, $\tilde{f} \circ \Gamma:I \to \Z$ is such a lift. By Proposition \ref{prop:Uniqueness of liftings of coverings}, $\tilde{f} \circ \Gamma = V$, so \[ \partial_0 f_*[\gamma] =  [V(1)] = [\tilde{f} \circ \Gamma(1)] = \tilde{f}_* \partial_0[\gamma],\] which proves the claim.
	
	The result now follows from the five-lemma. 
	\end{proof}

Using a homotopy equivalence directly, we have, for $r \geq \frac{1}{2}$,
\begin{theorem}
\label{thm:Fundamental group of circle large r}
$\pi_1(S^1,c_r) = \{1\}$ for $r \geq \frac{1}{2}$
\end{theorem} 
\begin{proof}
If $r \geq \frac{1}{2}$, then the closure structure $c_r$ is indiscrete, i.e. $c_r(A) = S^1$ for any $A \subset S^1$. Any function to $(S^1,c_r)$ is continuous in this case, and therefore $(S^1,c_r)$ is contractible. 
\end{proof}

Similarly,

\begin{theorem}
For $m \geq \left\lfloor\frac{n}{2}\right\rfloor$, or $n=3, m\geq 1$, $\pi_1(Z_n,c_m) = \{1\}$
\end{theorem}
\begin{proof}
	As in the previous theorem, in this regime the space $(Z_n,c_m)$ is indiscrete, and therefore contractible. The result follows.
\end{proof}

Finally, the following calculation of the fundamental group of the wedge of
circles with different closure structures follows easily from Theorems
\ref{thm:Group VK}, \ref{thm:Fundamental group of the circle 1}, and \ref{thm:Fundamental group of circle large r}.
\begin{theorem}
	Let $X = S^1_{r_1} \vee S^1_{r_2}$ be the wedge sum of two circles $[0,1]/0\sim 1$, where $r_1, r_2 \geq 0$. We endow each $S^1_{r_i}$ with the \v Cech closure operator $c_{r_i}$, $i\in \{1,2\}$, and for a subset $U \subset X$, we define a \v Cech closure operator $c_X(U) \coloneqq c_{r_1}(U \cap S^1_{r_1}) \cup c_{r_2}(U\cap S^1_{r_2})$. Then \begin{equation*}
	\pi_1(X) = \begin{cases} F_{a,b} & r_1, r_2< \frac{1}{3}\\
	\Z & r_i < \frac{1}{3}, \, r_j \geq \frac{1}{2}, \,i\neq j\\
	1 & r_1, r_2 \geq \frac{1}{2}
	\end{cases}
	\end{equation*} 
Where $F_{a,b}$ denotes the free group on two generators.
\end{theorem}
\begin{proof}
First suppose that $r_1,r_2 < \frac{1}{3}$. Recall that $S^1_{r_1} \vee
S^1_{r_2} = S^1_{r_1} \sqcup S^1_{r_2} / 0_1 \sim 0_2$ with the quotient
closure structure. (By Proposition \ref{prop:Quotient closure formula}, the
quotient closure structure agrees with the one given in the statement of the
theorem.) Define $U_1 = \langle S^1_{r_1} \sqcup (-r_2 -\epsilon, r_2+
\epsilon)\rangle$, for some $\epsilon >0$ small, where $\langle \cdot \rangle$
indicates the equivalence classes and all quantities are understood mod $1$.
Similarly, define $U_2 = \langle (-r_1 - \epsilon', r_1 + \epsilon') \sqcup
S^1_{r_2}$. Then $c_X(X-U_1) = [\epsilon, 1-\epsilon] \subset S^1_{r_2}$, and
therefore $S^1_{r_1} \subset \Int(U_1) = X - c_X(X-U_1)$. Similarly, $S^1_{r_2}
\subset \Int(U_2)$. Therefore, the collection $\mathcal{U} = \{U_1,U_2\}$ is an interior cover, and from the definitions, we see that the intersection $U_1 \cap U_2$ is contractible (and in particular) path-connected, and we apply Theorem \ref{thm:Group VK}.

For $r_i\geq \frac{1}{2}$, $r_j < \frac{1}{3}$, we use the same cover as above and apply Theorem \ref{thm:Group VK}. However, since $U_i$ is now contractible, $\pi_1(X) = \pi_1(U_2) = \Z$.

For $r_i \geq \frac{1}{2}$, $i = 1,2$, the closure structure $c_X$ is indiscrete. Therefore any map to $X$ is continuous, and $X$ is contractible. 
\end{proof}

\subsection{Combinatorial homotopy on graphs and simplicial complexes}
\label{subsec:HomSimp} 

In this section, we briefly show how to put a closure structure on graphs and, more generally, on $k$-skeleta of simplicial complexes.

Suppose $X$ is an abstract simplicial complex. Denote by $X^{(k)}$ the $k$-skeleton of $X$, $S(X)$ the set of simplices of $X$, and $S^{(k)}(X)$ the set of $k$-simplices of $X$. Let $\mathcal{P}(S^{(k)}(X)$ be the powerset of the set of $k$-simplices of $X$. We define an operator $c_{X,k}: \mathcal{P}(S^{(k)}(X)) \to \mathcal{P}(S^{(k)}(X))$ by 
\begin{equation}
\label{eq:Simplicial complex convergence structure}
c_{X,k}(\sigma) = \sigma \cup \{\gamma \in S^{(k)}(X) \mid \gamma \cap \sigma \in S^{(k-1)}(X) \text{ or } \gamma \cup \sigma \in S^{(k+1)}(X) \},
\end{equation} 
for $\sigma \in S^{(k)}(X)$. Now let
\begin{equation*}
c_{X,k}(A) \coloneqq \bigcup_{\sigma \in A} c_{X,k}(\sigma),
\end{equation*}
for $A \in \mathcal{P}(S(X))$.
Similarly, for the vertices $V$ in a directed graph $G = (V,E)$ we may define
\begin{equation}
c_{E}(U) = \{v \in V \mid v\in U \text{ or } (w,v)\in E \text{ for some } w\in U\}
\end{equation}

 We then have
\begin{proposition}
$c_{X,k}$ and $c_E$ are \v {C}ech closure operators.
\end{proposition}
\begin{proof} Immediate from the definitions of \v {C}ech closure operators, $c_{X,k}$, and $c_E$.
\end{proof}

For `circular graphs' on three vertices, we have the following.
\begin{theorem}
For the graph $G = (V,E)$ with vertices $V=\{0,1,2\}$, and edges $\{(k,k+1) \mod 3, k\in V\}$, the space $(V,c_E)$ is contractible.
\end{theorem}
\begin{proof}
	The closure structure is indiscrete, which implies that any map to $V$ is continuous. Therefore $V$ is contractible.
\end{proof}

As the spaces $(\Z_n,c_{n,m})$ can be interpreted as graphs by defining $V = \Z_n$, $E \coloneqq \{ (i,j) \mid j \in c(i)$, the calculations in Subsection \ref{subsec:Covering spaces} give a number of calculations of fundamental groups for graphs with the induced closure structure.

\begin{rem}Applying the above homotopy theory to the closure structures given
	by the $c_{X,k}$ for a simplicial complex $X$ results in homotopy
	groups in a similar spirit to those found in \cites{Babson_etal_2006,
	Barcelo_Capraro_White_2014, Barcelo_Laubenbacher_2005,
Barcelo_Kramer_Laubenbacher_2001}. Additionally, the homotopy theory for
$(V,c_E)$ closely resembles that given in \cite{Grigoryan_et_al_2014}. It is
not yet clear whether the homotopy groups developed here are the isomorphic to
those in \cites{Babson_etal_2006, Barcelo_Capraro_White_2014, Barcelo_Laubenbacher_2005, Barcelo_Kramer_Laubenbacher_2001} for graphs and simplicial complexes, or to those in \cite{Grigoryan_et_al_2014} for directed graphs.
\end{rem}

\subsection{Persistent homotopy}

We remark here that the closure structures given in Section \ref{sec:ScMet} can be used to define so-called `persistent' homotopy groups in the following way. If $q \leq r$, then the identity maps 
\begin{equation}
\label{eq:Iota}
\iota^{(q,r)}:(X,c_q) \to (X,c_r)
\end{equation} are continuous, and therefore induce maps $\iota^{(q,r)}_n:\pi_n(X,c_q) \to \pi_n(X,c_r)$, $n\geq 0$, which are homomorphisms for $n \geq 1$. We define the \emph{persistent homotopy groups} of a metric space $X$ by 
\begin{equation*}
\pi_n^{(q,r)}(X,d) = im \left(\iota_n^{(q,r)} \right), n\geq 1
\end{equation*} 
and the \emph{persistent components} of $X$ by \[\pi_0^{(q,r)}(X,d) = im \left(\iota_0^{(q,r)}\right).\]

\begin{rem}
We additionally note that, for any metric space $(X,d_X)$ and the collection of
closure structures $\{(X,c_r) \mid r \geq 0\}$, we may form a superlinear
family in the sense given in \cite{Bubenik_deSilva_Scott_2015} in the following
way. First, let $P=[0,\infty)$ be endowed with the standard order, and let
$\catname{D}=\catname{ScMet}$. The maps $\Omega_\epsilon: t \mapsto t+\epsilon$
is a superlinear family in $(P,\leq)$, and therefore, by Theorem 3.21 in
\cite{Bubenik_deSilva_Scott_2015}, the interleaving distance (Definition 3.20
in \cite{Bubenik_deSilva_Scott_2015}) gives a pseudometric $d^{\Omega}$ on the category $\catname{D^P}$ of functors from $\catname{P} \to \catname{ScMet}$. We are grateful to Peter Bubenik for this observation.
\end{rem}

\acknowledgement{We would like thank  H\'{e}l\`{e}ne Barcelo for interesting and useful discussions on her work on discrete homotopy, Conrad Plaut for discussions on his work in \cite{Plaut_Wilkins_2013}, Peter Bubenik for several useful remarks and for the reference to \cite{Bubenik_deSilva_Scott_2015}, Frank Weilandt for helpful comments on an earlier version of this paper, and Marian Mrozek and the Jagiellonian University for excellent hospitality during a short visit. We would also like to express our gratitude to the reviewer for the many detailed, thoughtful remarks which greatly improved the paper and helped to streamline a number of proofs.}

\bibliography{/home/antonio/Bib/all}
\end{document}